# STATISTICAL PERFORMANCE OF SUPPORT VECTOR MACHINES

By Gilles Blanchard,[1] Olivier Bousquet and Pascal Massart

*Fraunhofer-Insitute FIRST, Google and Université Paris-Sud*

The support vector machine (SVM) algorithm is well known to the computer learning community for its very good practical results. The goal of the present paper is to study this algorithm from a statistical perspective, using tools of concentration theory and empirical processes.

Our main result builds on the observation made by other authors that the SVM can be viewed as a statistical regularization procedure. From this point of view, it can also be interpreted as a model selection principle using a penalized criterion. It is then possible to adapt general methods related to model selection in this framework to study two important points: (1) what is the minimum penalty and how does it compare to the penalty actually used in the SVM algorithm; (2) is it possible to obtain "oracle inequalities" in that setting, for the specific loss function used in the SVM algorithm? We show that the answer to the latter question is positive and provides relevant insight to the former. Our result shows that it is possible to obtain fast rates of convergence for SVMs.

**1. Introduction.** The success of the support vector machine (SVM) algorithm for pattern recognition is probably mainly due to the number of remarkable experimental results that have been obtained in very diverse domains of application. The algorithm itself can be written as a nice convex optimization problem for which there exists a unique optimum, except in rare degenerate cases. It can also be expressed as the minimization of a regularized functional where the regularizer is the squared norm in a Hilbert space of functions on the input space. Although these are nice mathematical formulations, quite amenable to analysis, the statistical behavior of this

Received October 2006; revised April 2007.

[1]Supported in part by a grant of the Humboldt Foundation.

algorithm remains only partially understood. Our goal in this work is to investigate the properties of the SVM algorithm in a statistical setting.

1.1. *The abstract classification problem and convex loss approximation.* We consider a generic (binary) classification problem, defined by the following setting: assume that the product $\mathcal{X} \times \mathcal{Y}$ is a measurable space endowed with an unknown probability measure $P$, where $\mathcal{Y} = \{-1, 1\}$ and $\mathcal{X}$ is called the input space. The pair $(X, Y)$ denotes a random variable with values in $\mathcal{X} \times \mathcal{Y}$ distributed according to $P$. We will denote $P_X$ the marginal distribution of variable $X$. We observe a set of $n$ independent and identically distributed (i.i.d.) pairs $(X_i, Y_i)_{i=1}^n$ sampled according to $P$. These random variables form the *training set.*

Given this sample, the goal of the classification task is to estimate the *Bayes classifier*, that is, the measurable function $s^*$ from $\mathcal{X}$ to $\mathcal{Y}$ which minimizes the probability of misclassification, also called *generalization error*, $\mathcal{E}(s^*) = \mathbb{P}[s^*(X) \neq Y]$. It is easily shown that $s^*(x) = 2 \times \mathbf{1}\{P(Y = 1 | X = x) > \frac{1}{2}\} - 1$ a.s. on the set $\{P(Y = 1 | X = x) \neq \frac{1}{2}\}$. Note that it is an abuse to call $s^*$ "the" minimizer of the misclassification error, since it can have arbitrary value on the set $\{P(Y = 1 | X = x) = \frac{1}{2}\}$. In the sequel, we refer to $s^*$ as a fixed function, for example, if we choose arbitrarily $s^*$ to be 1 on the latter set.

Having a finite sample from $P$, a seemingly reasonable procedure is to find a classifier $s$ minimizing the empirical classification error $\mathcal{E}_n(s) = \frac{1}{n} \sum_i \mathbf{1}\{s(X_i) \neq Y_i\}$, with the minimization performed over some model of controlled complexity. However, this is in most cases intractable in practice because it is not a convex optimization procedure. This is the reason why a number of actual classification algorithms replace this loss by a convex loss over some real-valued (instead of $\{-1, 1\}$ valued) function spaces. This is the case of the SVM where such a "proxy" loss is used ensuring convexity properties. Its relation with the classification loss will be detailed in Section 2.1.

1.2. *Motivations.*

*Relative loss and oracle-type inequalities.* In the last two decades of the last century, the theoretical study of various classification algorithms has mainly focused on deriving confidence intervals about their generalization error. The foundations of this theory have been laid down by Vapnik and Chervonenkis as soon as 1971 [38]. Such confidence intervals have been derived for SVMs and, more generally, so-called "large margin classifiers," for example, using the notion of fat-shattering VC dimension; see [2].

However, it is probably fair to say that the explicit confidence intervals thus obtained are *never* sharp enough to be of practical interest—even though effort, legitimately, has been and is still made to obtain tighter bounds. On the other hand, we argue that uniform confidence intervals about the generalization error are not the most adapted tool to understand correctly the behavior of the algorithm.

If we compare the classification setting to regression, we see that, in regression, the loss of an estimator is always measured relatively to a target function $f^*$ (e.g., through $L^2$ distance). Furthermore, recent work (see, e.g., [22]) has shown that a precise study of the behavior of the relative loss when the estimator $\hat{f}$ is close to $f^*$ is a key element for proving correct convergence rates. This approach is sometimes called "localization."

In this paper we follow this general principle in the context of SVMs. Our main quantity of interest will therefore be the relative loss, for the proxy loss function, of $\hat{s}$ with respect to $s^*$, instead of the absolute loss itself (the average relative loss will also be called *risk*). In this regard, this work should be put in the context of a general trend in the recent literature on classification and, more generally, statistical learning, where the focus has shifted to the relative loss (see also below Section 5.1.2 for further discussion on this point).

Of course, a confidence interval for the relative loss is not informative, since $s^*$ is unknown; instead, the goal to be aimed at is an oracle-type inequality. The term *oracle inequality* originally refers to a risk bound for a model selection procedure where the bound is within a constant factor of the risk of a minimax estimator in the best model; that is, almost as good as if this best model had been known in advance through an "oracle". In the present context, we use more loosely the term "oracle-type inequality" to designate a bound where the risk of the estimator can be compared to the risk of the best approximating functions coming from any model under consideration plus a model-dependent penalty term; this without knowing in advance which models are best. This approach typically allows us to obtain precise bounds on the rates of convergence toward the target function.

*SVM and regularized model selection.* It has been noted by several authors (see Section 2.3) that SVMs can be seen as a regularized estimation method, where the regularizer is the squared norm of the estimating function in some reproducing kernel Hilbert space. We show that this can also be interpreted as a penalized model selection method, where the models are balls in this Hilbert space. This allows us to cast the SVM problem into a general penalized model selection framework, where we are able to use tools developed in [22], in order to obtain oracle-type inequalities over the family of considered models.

### 1.3. *Highlights of the present work.*

*A generic, versatile model selection theorem.* To be applied to SVMs, the results of [22] need to be extended to a setting where various parameters are model-dependent, resulting in various technical problems. Therefore, we decided to devote a whole section (Section 4) of this paper to the extension of these model selection results in a very general setting. We believe this result is of much interest per se because it can be useful for other applications (at least when the loss function is bounded model-wise) and constitutes an important point of this work.

*Is the SVM an adaptive procedure?*   The application of the above general result to SVMs is an example of the power of this approach, and allows us to derive a nonasymptotic oracle-type inequality for the SVM proxy risk. This is the main result of this paper. The interesting feature of oracle-type bounds is that they display adaptivity properties: while the regularization term used in the estimator does not depend on assumptions on the target function, the bound itself involves the approximation properties of the models to the target function. Therefore, the (fixed) estimation procedure "adapts" to how well the target is approximated by the models. This is in contrast to other related work on the subject such as [12, 32], where typically the optimal bound is obtained for a choice of the regularization constant that depends a priori on these approximation properties.

*Is the SVM regularization function adequate?*   Our result allows us to cast a new light on a very interesting problem, namely, concerning the adequate regularization function to be used in the SVM setting. Our main theorem establishes that the oracle-type inequality holds provided the regularizer function is larger than some lower bound $\zeta(\|f\|_k, n)$, which is a function of the Hilbert norm $\|f\|_k$ and the sample size $n$. Since the oracle inequality bound is nonincreasing in function of the regularization term, choosing the regularization precisely equal to $\zeta(\|f\|_k, n)$ will result in the best possible bound allowed by our analysis. The precise behavior, as a function of the sample size $n$, of $\zeta(\|f\|_k, n)$ depends on a capacity analysis of the kernel Hilbert space. For this, we provide two possible routes, either using the spectrum of the kernel integral operator, or the supremum norm entropy of the kernel space. In particular, we show (in both situations) that, while the *squared* Hilbert norm is traditionally used as a regularizer for the SVM, a *linear* function of the Hilbert norm is enough to ensure the oracle inequality: this suggests that the traditional regularizer could indeed be *too heavy*.

*Using several kernels.*   Another interesting consequence of the model selection approach is that it is possible to derive almost transparently an oracle-type inequality in an extended situation where we use several kernels at once for the SVM. Namely, the different kernels can be compared via their respective penalized empirical losses. The oracle inequality then states that this amounts to selecting the best kernel available for the problem.

*Influence of the generating probability on the convergence rate.* It has been recently pointed out (see [23, 35]) that in the classification setting, the behavior of the function $\eta(x) = \mathbb{P}[Y = 1 | X = x]$ in the neighborhood of the value $\frac{1}{2}$ plays a crucial role in the optimal convergence rate toward the Bayes classifier. In this paper we assume that $\eta(x)$ is bounded away from the value $\frac{1}{2}$ by a "gap" $\eta_0$ and study the influence of $\eta_0$ on the risk bounds obtained. An interesting feature of the result is that the knowledge of $\eta_0$ is not needed to define the estimator itself: it only comes into play through a remainder term in the bound.

Note that, for a strictly convex proxy loss, this type of assumption on $\eta$ essentially influences the relation between classification risk and proxy risk (see [4]), while it has no impact on the statistical behavior of the proxy risk itself. Because the proxy loss used by the SVM is not strictly convex (it is piecewise linear), the setting considered in the present paper is different: the gap assumption plays a role directly in the inequalities for the proxy risk and not in the relation with the classification risk.

1.4. *Organization of the paper.* In Section 2 we present the SVM algorithm, show how to formulate it as a model selection via penalization method and survey existing results. In Section 3 we state the main result of the paper for the SVM and discuss its implications and scope. The main tool to derive these results, which handle penalized model selection in a generic setting, is given in Section 4—we hope that its generality will make it useful in the future for other settings as well. We subsequently show how to apply this general result to the special case of the SVM. Section 5 contains a comparison of our result to other related work and concluding remarks. Finally, Section 6 contains the proofs of the results.

**2. Support vector machines.** For details about the algorithm, its basic properties and various extensions, we refer to the books [13, 29, 37]. We give here a short presentation of the formulation of the algorithm with emphasis on the fact that it can be thought of as a model selection via penalization method.

2.1. *Preliminaries*: *loss functions.* With some abuse of notation, we denote $Pg := \mathbb{E}[g(X, Y)]$ for an integrable function $g$ from $\mathcal{X} \times \mathcal{Y}$ to $\mathbb{R}$. Also, we introduce the empirical measure defined by the sample as

$$P_n := \frac{1}{n} \sum_{i=1}^{n} \delta_{X_i} \otimes \delta_{Y_i},$$

so that $P_n g$ denotes $n^{-1} \sum_{i=1}^{n} g(X_i, Y_i)$. Finally, we denote $\eta(x) = P[Y = 1 | X = x]$.

Before we delve further into the details of the support vector machine, we want to establish a few general preliminaries useful to understand the goals of the rest of the paper.

The natural setting to study SVMs is *real-valued classification* where we build estimators $\widehat{f}_n$ of $s^*$ as real-valued functions, being understood that the actual binary classifier associated to a real function is obtained by taking its sign. We therefore measure the probability of misclassification by comparing the sign of $\widehat{f}_n(X)$ to $Y$, thus rewriting the generalization error as

$$\mathcal{E}(\widehat{f}_n) = \mathbb{P}[Y\widehat{f}_n(X) \le 0] = \mathbb{E}[\theta(Y\widehat{f}_n(X))],$$

where $\theta(z) = \mathbf{1}\{z \le 0\}$ is called the 0-1 loss function. By a slight abuse of notation, we also denote by $\theta$ the following functional:

$$\theta(f) := (x,y) \mapsto \mathbf{1}\{yf(x) \le 0\}.$$

We define the associated risk (or relative average loss) function

$$\Theta(\widehat{f}_n, s^*) := \mathbb{P}[Y\widehat{f}_n(X) \le 0] - \mathbb{P}[Ys^*(X) \le 0] = P\theta(\widehat{f}_n) - P\theta(s^*).$$

However, as will appear in the next section, the classification error $\theta(\cdot)$ is not the actual measure of fit used by the algorithm of the support vector machine; it uses instead the "hinge loss" function defined by $\ell(z) := (1-z)_+$, where $(\cdot)_+$ denotes the positive part. Similarly, we also denote by $\ell$ the following functional:

$$\ell(f) := (x,y) \mapsto (1 - yf(x))_+;$$

the associated risk function is denoted

$$L(\widehat{f}_n, s^*) := \mathbb{E}[\ell(\widehat{f}_n) - \ell(s^*)].$$

As mentioned in the introduction, using this convex loss allows for a tractable optimization problem for actual implementation of the algorithm. Since $\ell$ is the loss function actually used to build the SVM classifier, the aim of our analysis is to derive oracle inequalities about its associated risk $L$.

However, as the main goal of classification is ultimately to obtain low generalization error $\mathcal{E}$, it is only natural to ask the question of the connection between the two above losses. It is obvious that $\theta(x) \le \ell(x)$ and therefore that $\mathcal{E}(f) \le \mathbb{E}[\ell(f)]$. Nevertheless, recalling our main focus is on risks (i.e., relative average loss), this remark is not really satisfactory and the two following additional questions are of primary interest:

- How is the *real-valued* function $f^*$ minimizing the averaged hinge loss $\mathbb{E}[\ell(f^*)]$ related to the optimal classifier $s^*$?
- How are $\Theta(\cdot, \cdot)$ and $L(\cdot, \cdot)$ related?

(Again, note that it is not entirely correct to talk about "the" function $f^*$ minimizing the hinge loss, since it is not unique: in the sequel we will assume a specific choice has been made.)

The following elementary lemma gives a satisfactory answer to these questions:

LEMMA 2.1. (i) *Let $s^*$ be a minimizer of $\mathcal{E}(s)$ over all measurable functions $s$ from $\mathcal{X}$ into $\{-1, 1\}$. Then the following holds:*

$$\mathbb{E}[\ell(s^*)] = \min_f \mathbb{E}[\ell(f)],$$

*where the right-hand side minimum is taken over all measurable real-valued functions on $\mathcal{X}$. Furthermore, if $f^*$ is a minimizer of $E[\ell(f)]$, then $f^* = s^*$ a.s. on the set $\{\mathbb{P}[Y = 1 | X = x] \notin \{0, \frac{1}{2}, 1\}\}$.*

(ii) *For any $P$-measurable function $f$,*

$$\Theta(f, s^*) \leq L(f, f^*).$$

Part (i) of the lemma can be found in [19] and part (ii) in [40], but we give a self-contained proof in Section 6.1 for completeness. Since the choice of $f^*$ is arbitrary among minimizers of $\mathbb{E}[\ell(f)]$, (i) implies that we can choose $f^* = s^*$, which will be assumed from now on.

2.2. *The SVM algorithm.* There are several possible ways of formulating the SVM algorithm. Historically, it was formulated geometrically. First suppose the input space $\mathcal{X}$ is a Hilbert vector space and that the two classes can be separated by a hyperplane. The SVM classifier is then the linear classifier obtained by finding the hyperplane which separates the training points in the two classes with the largest margin (maximal margin hyperplane). The margin corresponds to the smallest distance from a data point to the hyperplane.

Now, in general, $\mathcal{X}$ may not be a Hilbert space, but is *mapped* into one where the above algorithm is applied. For computational tractability of the algorithm, it is crucial that this Hilbert space can be generated by a (reproducing) kernel, whose properties we sum up briefly here.

Assume we have at hand a so-called kernel function $k: \mathcal{X} \times \mathcal{X} \to \mathbb{R}$, meaning that $k$ is symmetric and positive semi-definite, in the following sense:

$$\forall n, \forall (x_1, \ldots, x_n) \in \mathcal{X}^n, \forall (a_1, \ldots, a_n) \in \mathbb{R}^n \qquad \sum_{i,j=1}^n a_i a_j k(x_i, x_j) \geq 0.$$

It can be proved that such a function defines a unique reproducing kernel Hilbert space (RKHS for short) $\mathcal{H}_k$ of real-valued functions on $\mathcal{X}$. Namely, define $\mathcal{H}_k$ as the completion of $\mathrm{span}\{k(x,\cdot):x\in\mathcal{X}\}$, with respect to the norm induced by the following inner product:

$$\langle u,v\rangle_k=\sum_{i=1}^{n}\sum_{j=1}^{m}a_i b_j k(x_i,x_j)\qquad\text{for }u=\sum_{i=1}^{n}a_i k(x_i,\cdot)\text{ and }v=\sum_{j=1}^{m}b_j k(x_j,\cdot);$$

here the completion is defined in such a way so that it consists of real functions on $\mathcal{X}$ as announced. We denote the norm in $\mathcal{H}_k$ by $\|\cdot\|_k$.

Since $\mathcal{H}_k$ is a Hilbert space of real-valued functions on $\mathcal{X}$, any element $w$ of $\mathcal{H}_k$ can be alternatively understood as a vector or as a function. Moreover, this space has the so-called reproducing property which can be expressed as

$$\forall u\in\mathcal{H}_k,\forall x\in\mathcal{X}\qquad u(x)=\langle u,k(x,\cdot)\rangle_k.$$

Finally, as announced, the input space $\mathcal{X}$ is mapped into $\mathcal{H}_k$ by the simple mapping $x\mapsto k(x,\cdot)$, and, thus, the scalar product of the images of $x,x'\in\mathcal{X}$ in $\mathcal{H}_k$ is just given by $k(x,x')$.

Now, in that space, a hyperplane is defined by its normal vector $w$ and a threshold $b\in\mathbb{R}$ as

$$H(w,b)=\{v\in\mathcal{H}_k:\langle w,v\rangle_k+b=0\}.$$

It is easy to see [29] that the maximum margin hyperplane (when it exists) is given by the solution of the following optimization problem:

$$\min_{w\in\mathcal{H}_k,b\in\mathbb{R}}\tfrac{1}{2}\|w\|_k^2$$

under the constraints: $\forall i=1,\dots,n,Y_i(\langle w,k(X_i,\cdot)\rangle_k+b)\geq 1.$

However, it can happen that the data is not linearly separable (i.e., the above constraints define an empty set). This has led to considering the following relaxed optimization problem, depending on some constant $C\geq 0$:

$$\min_{w\in\mathcal{H}_k,b\in\mathbb{R}}\tfrac{1}{2}\|w\|_k^2+C\sum_{i=1}^{n}\xi_i$$

(2.1) under the constraints: $\forall i=1,\dots,n,Y_i(\langle w,k(X_i,\cdot)\rangle_k+b)\geq 1-\xi_i;$

$$\forall i=1,\dots,n\qquad\xi_i\geq 0.$$

This problem always has a solution and is usually referred to as the soft-margin SVM. It is common, although not systematical, for theoretical studies of SVMs to introduce a simpler version of the SVM algorithm where one uses only hyperplanes containing the origin, that is, $b$ is set to zero (although this version is admittedly rarely used in practice). This is mainly for avoiding some technical difficulties. We will adopt this simplification here, calling this constrained version "SVM$_0$," and we will focus on it for the main part of the paper.

2.3. *From regularization to model selection.* It has been noticed by several authors [15, 30] that the soft-margin SVM algorithm can be formulated as the minimization of a regularized functional. Consider the primal optimization problem (2.1). For a fixed $w$, obviously the optimal choice for the parameters $(\xi_i)$ given the constraints is $\xi_i = (1 - Y_i(\langle w, k(X_i, \cdot) \rangle_k + b))_+$. Now using the reproducing property of the kernel, we have $\langle w, k(X_i, \cdot) \rangle_k = w(X_i)$, so the new formulation of the problem is (now denoting $f$ instead of $w$)

$$(2.2) \qquad \min_f \frac{1}{n} \sum_{i=1}^n (1 - Y_i f(X_i))_+ + \Lambda_n \|f\|_k^2,$$

where $\Lambda_n = \frac{1}{nC}$ and the minimum is to be performed over $f \in \mathcal{H}_k$ (for the SVM$_0$ algorithm) and for $f \in \mathcal{H}_k^b = \{x \mapsto g(x) + b | g \in \mathcal{H}_k, b \in \mathbb{R}\}$ for the plain SVM algorithm. Note that $\|\cdot\|_k$, inherited from $\mathcal{H}_k$ to $\mathcal{H}_k^b$, is only a semi-norm on $\mathcal{H}_k^b$.

Now, it is straightforward that the optimization problem (2.2) can be rewritten in the following way:

$$(2.3) \qquad \min_{R \in \mathbb{R}} \left\{ \min_{f : \|f\|_k \le R} \frac{1}{n} \sum_{i=1}^n (1 - Y_i f(X_i))_+ + \Lambda_n R^2 \right\}.$$

This gives rise to the interpretation of the above regularization as *model selection*, where the models are balls in $\mathcal{H}_k$ (or "semi-norm balls" in $\mathcal{H}_k^b$), and where the model selection is done using *penalized empirical loss minimization.* Also, it is now clear from equations (2.2) and (2.3) that the empirical loss used by the SVM is not the classification error (or 0–1 loss function), but the hinge loss function $\ell$ defined in the previous section.

Denoting $\mathcal{B}(R)$ the ball of $\mathcal{H}_k$ of radius $R$, our interest in the main part of the paper is to study the behavior of SVM$_0$ *vis-à-vis* the family of models $\mathcal{B}(R)$, and the correct order of the regularization function to be used.

## 3. Main result.

3.1. *Assumptions.* We will present two variations of our main result. The difference between the two versions is in the way the capacity of the RKHS is analyzed. General assumptions on the RKHS $\mathcal{H}_k$ and on the generating distribution are common to the two versions. Below we denote $\eta(x) = P(Y = 1 | X = x)$.

(A1) $\mathcal{H}_k$ is a separable space (Note that the separability of $\mathcal{H}_k$ is ensured, in particular, if $\mathcal{X}$ is a compact topological space and $k$ is continuous on $\mathcal{X} \times \mathcal{X}$.), and $k(x, x) \le M^2 < \infty$ for all $x \in \mathcal{X}$.

(A2) ("Low noise" condition) $\qquad \forall x \in \mathcal{X} \qquad |\eta(x) - \frac{1}{2}| \ge \eta_0.$

The following additional assumption will be required only for setting (S1) below:

(A3)                    $\forall x \in \mathcal{X}$        $\min(\eta(x), 1 - \eta(x)) \geq \eta_1$.

Our result covers the two following possible settings:

*Setting* 1 (S1).  Suppose assumptions (A1), (A2) and (A3) satisfied. In this first setting, the capacity of the RKHS is analyzed through the spectral properties of the kernel integral operator $L_k : L^2(P_X) \to L^2(P_X)$ defined as

(3.1)                    $$(L_k f)(x) = \int k(x, x') f(x') \, dP_X(x'),$$

which is positive, self-adjoint and trace-class (see Appendix A for details). As a result, $L_k$ can be diagonalized in an orthogonal basis of $L^2(P_X)$, it has discrete spectrum $\lambda_1 \geq \lambda_2 \geq \cdots$ (where the eigenvalues are repeated with their multiplicities) and satisfies $\sum_{j \geq 0} \lambda_j < \infty$. For a fixed $\delta > 0$, we then define for $n \in \mathbb{N}$ the following function:

$$\gamma(n) = \eta_1^{-1} \frac{1}{\sqrt{n}} \inf_{d \in \mathbb{N}} \left( \frac{d}{\sqrt{n}} + \frac{\eta_1}{M} \sqrt{\sum_{j > d} \lambda_j} \right).$$

*Setting* 2 (S2).  Suppose assumptions (A1) and (A2) satisfied. For the second situation covered by the theorem, the capacity is measured via supremum norm covering numbers. In this situation, we assume that the RKHS $\mathcal{H}_k$ can be included via a compact injection into $\mathcal{C}(\mathcal{X})$ and we denote by $H_\infty(\mathcal{B}_{\mathcal{H}_k}, \varepsilon)$ the $\varepsilon$-entropy number (log-covering number) in the supremum norm of the unit ball of $\mathcal{H}_k$. Denote

(3.2)                    $$\xi(x) = \int_0^x \sqrt{H_\infty(\mathcal{B}_{\mathcal{H}_k}, \varepsilon)} \, d\varepsilon,$$

and let $x_*(n)$ be the solution of the equation $\xi(x) = M^{-1} n^{1/2} x^2$. For a fixed $\delta > 0$, define for $n \in \mathbb{N}$ the following function:

$$\gamma(n) = M^{-2} x_*^2(n).$$

3.2. *Statement.*  We now state our main result, which applies, in particular, to the SVM$_0$ algorithm.

THEOREM 3.1.  *Consider either setting* (S1) *under assumptions* (A1), (A2) *and* (A3), *or setting* (S2) *under assumptions* (A1) *and* (A2). *Define the constant* $w_1 = \eta_1$ *for setting* (S1) *and* $w_1 = 1$ *for setting* (S2).

*Let* $\delta > 0$ *be a fixed real number; and let* $\Lambda_n > 0$ *be a real number satisfying*

(3.3)                    $$\Lambda_n \geq c \left( \gamma(n) + w_1^{-1} \frac{\log(\delta^{-1} \log n) \vee 1}{n} \right),$$

*where $c$ is a universal constant. Finally, let $\varphi$ be a nondecreasing function on $\mathbb{R}^+$ such that $\varphi(0) = 0$ and $\varphi(x) \geq x$ for $x \geq \frac{1}{2}$.*

*Consider the following regularized minimum empirical loss procedure on an i.i.d. sample $((X_i, Y_i))_{i=1,\ldots,n}$ from distribution $P$, using the hinge loss function $\ell(x, y) = (1 - xy)_+$:*

$$(3.4) \qquad \widehat{g} = \operatorname*{Arg\,Min}_{g \in \mathcal{H}_k} \left( \frac{1}{n} \sum_{i=1}^{n} \ell(g(X_i), Y_i) + \Lambda_n \varphi(M\|g\|_k) \right),$$

*then if $s^*$ denotes the Bayes classifier, the following bound holds with probability at least $1 - \delta$:*

$$(3.5) \quad L(\widehat{g}, s^*) \leq 2 \inf_{g \in H_k} \left[ L(g, s^*) + 2\Lambda_n \varphi(2M\|g\|_k) \right] + 4\Lambda_n (2\varphi(2) + c w_1 \eta_0^{-1}).$$

### 3.3. *Discussion and comments.*

#### 3.3.1. *Discussion of the result.*

*Adaptivity of the SVM.* The most important point we would like to stress about Theorem 3.1 is that the regularization term and the final bound are *independent* of any assumption on how well the target function $f^*$ is approximated by functions in $\mathcal{H}_k$. This is an important advantage in the approach we advocate here, that is, casting regularization as model selection. The model selection approach dictates a minimal order of the regularization, which is "structural" in the sense that it depends on some complexity measure of the models (here balls of $\mathcal{H}_k$) and not on how well the models approximate the target. In simpler terms, the minimal regularizer depends only on the estimation error, not the approximation error. Our result is therefore an *oracle type* bound, which entails that the SVM is an *adaptive* procedure with respect to the approximation properties of the target by functions in $\mathcal{H}_k$. From this bound, we can derive convergence rates to Bayes as soon as we have an additional hypothesis on these approximation properties, while the procedure stays unchanged. We discuss this point in more detail in Section 3.4.

*Squared versus linear regularization.* The second point we want to emphasize about Theorem 3.1 is that the minimum regularization function required to ensure that the oracle inequality holds is of order $\|g\|_k$ only (as a function of $\|g\|_k$). In the original SVM algorithm, a regularization of order $\|g\|_k^2$ is used. The theorem covers both situations by choosing respectively $\varphi(x) = x$ or $\varphi(x) = 2x^2$. In view of the oracle inequality, the weaker the regularization term, the better the upper bound: provided that the oracle inequality holds, a weaker regularization will grant a better bound on the convergence rate. Therefore, this theorem suggests that [under certain conditions, i.e., mainly (A2)] a lighter regularization can be used instead of the standard, quadratic, one.

Of course, while a lighter regularization results in a better *bound* in our theorem, we cannot assert positively that the resulting algorithm will necessarily outperform the standard one: to draw such a conclusion, we would need a corresponding lower bound for the standard algorithm. Here we will merely point out the analogy of SVM to regularized least squares regression. Under a Gaussian noise assumption, the behavior of the regularized least squares estimator of the form (3.4) [with the square loss $\ell(x, y) = (x - y)^2$ replacing the hinge loss] is completely elucidated (see [24], Section 4.4). In particular, the standard quadratic regularization estimator has an explicit form, from which it is relatively simple to derive corresponding lower bounds. As a consequence, in that case, it can be proven that a regularization that is lighter than quadratic enjoys better adaptivity properties than the standard one. In the present work, we have followed essentially the same driving ideas to derive our main result in the SVM setting, so that there is reasonable hope that the obtained bound indeed reflects the behavior of the algorithm. A complete proof of that fact is an interesting open issue.

*From hinge loss risk to classification risk.* This theorem relates the relative hinge loss $\mathbb{E}[\ell(\hat{g}) - \ell(s^*)]$ (where $s^*$ is the Bayes classifier) to the optimum relative loss in the models considered, that is, balls of $\mathcal{H}_k$ (see Section 2.3). Furthermore, Lemma 2.1 ensures that the relative classification error is upper-bounded by the relative hinge loss error, hence, the theorem also results in a bound on the relative classification error.

### 3.3.2. *Discussion of the assumptions.*

*About assumption* (A2). This assumption requires that the conditional probability of $Y$ given $X$ should be bounded away from $\frac{1}{2}$ by a "gap" $\eta_0$. Note that the knowledge of $\eta_0$ *is not necessary* for the definition of the estimator, as it does not enter in the regularization term. This quantity only appears as an additional term in the oracle inequality (3.5). Furthermore, for $\eta_0$ not depending on $n$, this trailing term will become negligible as $n \to \infty$, since the infimum in the first term will be attained for a function $g_n \in \mathcal{H}_k$ with $\|g_n\|_k \to \infty$ (see below Section 3.4 ). Assumption (A2) is a particular case of the so-called Tsybakov's noise condition, which is known to be a crucial factor for determining fast minmax rates in classification problems (see [23, 35]).

*A possible generalization.* A more general Tsybakov's noise condition would be to assume, in place of (A2), that $|\frac{1}{2} - \eta(x)|^{-1} \in L^p$ for some $p > 0$. In this setting, it is possible to show (although it is out of the scope if the present work) that a result similar to (3.5) holds, with the same regularization function, except that the trailing term in (3.5) of order $\eta_0^{-1} \Lambda_n$ gets replaced by a term of the form $\zeta(\Lambda_n)$, with $x \lesssim \zeta(x) \lesssim \sqrt{x}$, where the exact form of $\zeta$ depends on the noise condition and the structural complexity analysis of $\mathcal{H}_k$. Obviously, in this general situation, the trailing term is no longer

necessarily negligible — whether or not this is the case will depend on the behavior of the first term of the bound, and therefore on the approximation properties of $f^*$ by $\mathcal{H}_k$. The interpretation of this generalization is therefore more involved.

*About assumption* (A3). The requirement that $\eta$ should be bounded away from $0, 1$ by a gap $\eta_1$ is a technical assumption in setting (S2) needed as a *quid pro quo* for obtaining an explicit relation between regularization term and eigenvalues (see the short discussion before Theorem 6.6 in Section 6.3). While there does not appear to be an intrinsic reason for this assumption, we did not succeed in getting rid of it in this setting. Note that, in contrast to the previous point, the knowledge of $\eta_1$ is needed to define the regularization explicitly in this setting. While this assumption is somewhat unsatisfactory, it is possible, at least in principle, to obtain an explicit lower bound on the value of $\eta_1$ by introducing deliberately in the data a small artificial "label flipping noise" (i.e., flipping a small proportion of the training labels). We refer to [9] (in the discussion preceding Corollary 10 there; the idea also appeared earlier in [39]) where this idea is exposed in more detail. Note that the label flipping preserves assumption (A2), albeit with a smaller gap value $\eta_0$.

*About setting* (S2). An unsatisfactory part of the result for setting (S2) is that it is not possible to compute the value of the regularization parameter $\gamma(n)$ from the data, since it requires knowledge of the eigenvalues of $L_k$. The interest of this setting is to give an idea of what the relevant quantities are for defining a suitable regularization, in a way that is generally more precise than for setting (S1) (see discussion in the next section). Moreover, there is strong hope that estimating these tail sum of eigenvalues from the data (using, e.g., techniques from [3]) would lead to a suitable data-dependent penalty.

### 3.3.3. *Other comments.*

*Multiplicative constant.* The constant 2 in front of the right-hand side of equation (3.5) could be made arbitrarily close to one at the price of increasing the regularization function accordingly. Here we made an arbitrary choice in order to simplify the result.

*Deviation inequality vs. average risk.* The above result states a deviation bound valid with high probability $1 - \delta$. Note that $\delta$ enters into the regularization function, hence, it is not possible to directly integrate (3.5) to state a bound for the average risk. However, it is possible to obtain such a result at the price of a slightly heavier regularization (an additional logarithmic factor). Namely, the proof of Theorem 2 essentially relies on a general model selection theorem (Theorem 4.3 in the next section) which covers both the deviation inequalities and average risk inequalities with minor changes in the penalty function. For brevity, we do not state here the resulting theorem obtained for average SVM performance, but it should be clear that only minor modifications to the proof of Theorem 3.1 would be necessary.

*Using several kernels at once.* Suppose we have several different kernels $k_1, \ldots, k_t$ at hand. Then we can adapt the theorem to use them simultaneously. Namely, to each kernel is associated a penalization constant $\Lambda_n^{(i)}$; the estimator $\widehat{g}$ is given by (3.4) where we add another Arg min operation over the kernel index; and oracle inequality (3.5) is valid with an additional minimum over the kernel index; only $\delta$ has to be replaced by $\delta/t$ for the price of the union bound. That such a result holds is straightforward when one takes a look at the model selection approach used to prove Theorem 3.1 (developed Section 4). This is one of the advantages of this approach.

### 3.4. *Penalty functions and convergence rates for support vector machines.*

3.4.1. *Convergence rates for the SVM.* Let us first note from the definition of $\gamma$ in both settings (S1) and (S2) that, generally, $\gamma(n)$ is of order lower than $n^{-1/2}$. This is in contrast with some earlier results in learning theory where bounds and associated penalties often behave like $n^{-1/2}$. Actual rates of convergence to the Bayes classifier also depend on the behavior of the bias (or approximation error) term $\inf_{\|g\|_k \leq R} L(g, s^*)$. In most practical cases, the functions in $\mathcal{H}_k$ are continuous, while the Bayes classifier is not; hence, the Bayes classifier cannot belong to any of the models. If we assume that $\mathcal{H}_k$ is dense in $L_1(P)$, however (see also the stronger notion of "universal kernel" in [31]), then there exists a sequence of functions $(g_n) \in \mathcal{H}_k$ such that $u_n = L(g_n, s^*) \to 0$, implying consistency of the SVM. Moreover, if information is available about the speed of approximation [i.e., how $\inf_{\|g\|_k \leq R} L(g, s^*)$ goes to zero as a function of $R$] and about the function $\gamma(n)$ [depending either on eigenvalues or supremum norm entropy according to setting (S1) or (S2)], an upper bound on the speed of convergence of the estimator can be derived from Theorem 3.1. As noted earlier, in this case, using a regularization term of order $\|g\|_k$ instead of $\|g\|_k^2$ always leads to a better upper bound on the convergence rate. The study of such approximation rates for special function classes is outside the scope of the present paper, but is an interesting future direction.

3.4.2. *About the function $\gamma(n)$ in settings (S1) and (S2).* The behavior, as a function of the norm $\|g\|_k$, of the minimum regularization function required in the theorem does not depend on the setting. Its behavior as a function of the sample size $n$, however, does, since the complexity analysis is different in both settings.

In order to fix ideas, we give here a very classical Sobolev space type example where we can explicitly compute the function $\gamma$ in both settings— and where they coincide. Let us consider the case where $\mathcal{X} = \mathbb{T}$ is the unit circle, the marginal $P_X$ of the observations is the Lebesgue measure, and the reproducing kernel $k$ is translation invariant, $k(x, y) = k(x - y)$ where $k$ is a periodic function that admits the Fourier series decomposition

$$k(z) = \sum_{k \geq 0} a_k \cos(2\pi k z),$$

where $(a_k)$ is a sequence of nonnegative numbers. Obviously, the Fourier basis forms a basis of eigenvectors for the associated integral operator $L_k$ and the eigenvalues are $\lambda_1 = a_0, \lambda_{2k} = \lambda_{2k+1} = a_k/2$ for $k > 0$. A function belonging to the RKHS $f \in \mathcal{H}_k$ is therefore characterized by $\sum_{k \geq 0} \lambda_k^{-1} \widehat{f}_k^2 = \|f\|_k^2 < \infty$, where $\widehat{f}_k$ are its Fourier coefficients.

Consider the case where $\lambda_k \lesssim k^{-2s}$ for some $s > \frac{1}{2}$. Then computing the function $\gamma$ in setting (S1) yields $\gamma_1(n) \lesssim n^{-2s/(2s+1)}$. On the other hand, clearly $\mathcal{H}_k$ can be continuously included into the Sobolev space $H^s(\mathbb{T})$. Uniform norm entropy estimates for Sobolev spaces have been established (and can be traced back to [7]; see also [14], page 105 for a general result); it is known that $H_\infty(\mathcal{B}_{H^s(\mathbb{T})}, \varepsilon) \lesssim \varepsilon^{-1/s}$; hence, the function $\xi$ appearing in setting (S2) is such that $\xi(x) \lesssim n^{(2s-1)/2s}$, leading also to $\gamma_2(n) \lesssim n^{-2s/(2s+1)}$.

However, the fact that the two settings lead to a regularization of the same order seems very specific to this case, depending, in particular, on the properties of the Lebesgue measure and of Sobolev spaces. In a more general situation, if we assume the eigenvalues to be known, and $\eta_1$ to be a fixed constant, we expect the analysis in setting (S1) to give a tighter estimate for the minimal regularization function than the analysis in setting (S2); that is to say, the function $\gamma(n)$ appearing in (S1) will be of smaller order than the one appearing in (S2). Informally speaking, this is because the eigenvalues of $L_K$ are related to the covering entropy of the unit ball of $\mathcal{H}_k$ in $L^2$ norm, while setting (S2) considers covering entropy with respect to the stronger supremum norm.

On the other hand, this tighter analysis comes at a certain price, namely, additional assumption (A3) and the requirement that the eigenvalues are known (or estimated), as already pointed out above. One advantage of supremum norm entropy is that, by definition, it is distribution independent. Furthermore, some relatively general results are known on this entropy depending on the regularity properties of the kernel function; see [41].

## 4. A model selection theorem and its application.

### 4.1. *An abstract model selection theorem.* The remainder of the paper is devoted to the proof of Theorem 3.1. However, in the present section we change gears somewhat, forgetting voluntarily about the specific setting of

the SVM to present an "abstract" theorem resulting in oracle inequalities that can be obtained for model selection by penalized empirical loss minimization. This theorem is the cornerstone for the proof of Theorem 3.1.

Our motivation for leaving momentarily the SVM framework for a more general one is twofold. On the one hand, we hope that it will make appear more clearly to the reader the general principle underlying our result, independently of the specifics of the SVM (which we will return to in the next section). On the other hand, we think that this result is general enough to be of interest of itself, inasmuch as it can be applied in a variety of different frameworks.

The theorem is mainly an extended version of Theorem 4.2 of [22] to a more general setting, namely, where some key parameters, considered fixed in the above reference, can now depend on the model. This extension is necessary for our intended application to SVMs, which is exposed in Section 4.2, and requires appropriate handling. However, the scope of this abstract model selection theorem can cover a wider variety of situations. Examples are the classical VC-dimension setting using classification loss (in this case the result of [22] is actually sufficient; see also the more detailed study [23]), or regularized Boosting-type procedures (see [9], where an earlier version of the model selection theorem presented here was used). The fact that the theorem applies to approximate, rather than exact, penalized minimum empirical loss estimation is a minor refinement that is useful in certain situations: this will be the case for our application to SVMs, where the continuous regularization scheme will be related to an approximate discrete penalization scheme.

We first need to introduce the following definition:

DEFINITION 4.1. A function $\psi : [0, \infty) \to [0, \infty)$ is *sub-root* if it is nonnegative, nondecreasing, and if $r \mapsto \psi(r)/\sqrt{r}$ is nonincreasing for $r > 0$.

Sub-root functions have the following property:

LEMMA 4.2 ([3]). *Let* $\psi : [0, \infty) \to [0, \infty)$ *be a sub-root function. Then it is continuous on* $[0, \infty)$ *and the equation* $\psi(r) = r$ *has a unique positive solution. If we denote this solution by* $r^*$, *then for all* $r > 0$, $r \geq \psi(r)$ *if and only if* $r^* \leq r$.

We can now state the model selection result:

THEOREM 4.3. *Let* $\ell : \mathfrak{G} \to L^2(P)$ *[where* $\mathfrak{G} \subset L^2(P)$*] be a loss function and assume that there exists* $g^* \in \text{Arg\,Min}_{g \in \mathfrak{G}} \mathbb{E}[\ell(g)]$. *Let* $(\mathcal{G}_m)_{m \in \mathcal{M}}$, $\mathcal{G}_m \subset \mathfrak{G}$ *be a countable collection of classes of functions and assume there exists the following:*

- *a pseudo-distance $d$ on $\mathfrak{G}$;*
- *a sequence of sub-root functions $(\phi_m), m \in \mathcal{M}$;*
- *two positive sequences $(b_m)$ and $(C_m), m \in \mathcal{M}$;*

*such that*

(H1) $\qquad \forall m \in \mathcal{M}, \forall g \in \mathcal{G}_m \qquad \|\ell(g)\|_\infty \leq b_m;$

(H2) $\qquad \forall g, g' \in \mathfrak{G} \qquad \mathrm{Var}(\ell(g) - \ell(g')) \leq d^2(g, g');$

(H3) $\qquad \forall m \in \mathcal{M}, \forall g \in \mathcal{G}_m \qquad d^2(g, g^*) \leq C_m L(g, g^*);$

*and, if $r_m^*$ denotes the solution of $\phi_m(r) = r/C_m$,*

(H4) $\qquad \forall m \in \mathcal{M}, \forall g_0 \in \mathcal{G}_m, \forall r \geq r_m^*$

$$\mathbb{E}\left[\sup_{\substack{g \in \mathcal{G}_m \\ d^2(g,g_0) \leq r}} (P - P_n)(\ell(g) - \ell(g_0))\right] \leq \phi_m(r).$$

*Let $(x_m)_{m \in \mathcal{M}}$ be a sequence of real numbers such that $\sum_{m \in \mathcal{M}} e^{-x_m} \leq 1$. We assume that families $(b_m), (C_m), (x_m), m \in \mathcal{M}$, are ordered the same way, by which we mean that*

$$(4.1) \qquad \forall m, m' \in \mathcal{M}, \qquad x_m < x_{m'} \Rightarrow \begin{cases} b_m \leq b_{m'}; \\ C_m \leq C_{m'}. \end{cases}$$

*Let $\xi > 0, K > 1$ be some real numbers to be fixed in advance. Put $B_m = 75KC_m + 28b_m$, and let $\mathrm{pen}(m)$ be a penalty function such that, for each $m \in \mathcal{M}$,*

$$(4.2) \qquad \mathrm{pen}(m) \geq 250K \frac{r_m^*}{C_m} + \frac{B_m(x_m + \xi + \log(2))}{3n}.$$

*Let $(\rho_m)_{m \in \mathcal{M}}$ be a family of positive numbers and $\widetilde{g}$ denote a $(\rho_m)$-approximate penalized minimum empirical loss estimator over the family $(\mathcal{G}_m)$ using the above penalty function, that is, satisfying*

$$(4.3) \qquad \begin{aligned} &\exists \widetilde{m} \in \mathcal{M} : \widetilde{g} \in \mathcal{G}_{\widetilde{m}} \quad \text{and} \\ &P_n \ell(\widetilde{g}) + \mathrm{pen}(\widetilde{m}) \leq \inf_{m \in \mathcal{M}} \inf_{g \in \mathcal{G}_m} (P_n \ell(g) + \mathrm{pen}(m) + \rho_m); \end{aligned}$$

*then the following deviation inequality holds with probability greater than $1 - \exp(-\xi)$:*

$$L(\widetilde{g}, g^*) \leq \frac{K + 1/5}{K - 1} \inf_{m \in \mathcal{M}} \left( \inf_{g \in \mathcal{G}_m} L(g, g^*) + 2\mathrm{pen}(m) + \rho_m \right).$$

*Furthermore, if the penalty function satisfies, for each $m \in \mathcal{M}$,*

$$(4.4) \qquad \mathrm{pen}(m) \geq 250K \frac{r_m^*}{C_m} + \frac{B_m(x_m + \log(2))}{3n} + \frac{B_m \log B_m}{n},$$

*then the following expected risk inequality holds:*

$$\mathbb{E}[L(\widetilde{g}, g^*)] \leq \frac{K + 1/5}{K - 1} \inf_{m \in \mathcal{M}} \left( \inf_{g \in \mathcal{G}_m} L(g, g^*) + 2\mathrm{pen}(m) + \rho_m + \frac{2}{n} \right).$$

*Remarks.*

1. Note that the difference with Theorem 4.2 of [22] is the fact that constants $b_m$ and $C_m$ can depend on $m$, which requires additional work, but is a necessary step for application to SVMs.

2. In hypothesis (H4) $\phi(r^2)$ can be interpreted as the *modulus of continuity* with respect to $d$ of the supremum of the empirical process indexed by $\mathcal{G}$.

3. The class $\mathfrak{G} \subset L^2(P)$ should be seen as the "ambient space"; it should at least contain all models. Note that choice of $\mathfrak{G}$ determines the target function $g^*$ (the minimizer of the average loss on $\mathfrak{G}$). Typically, the theorem will be applied with $\mathfrak{G} = L^2(P)$ or $\mathfrak{G} = L^2(P_X)$ (as will be the case below), but other choices may be useful.

4. Although it is not its main purpose, this theorem can also be used for the convergence analysis of the empirical loss minimization procedure on a single model $\mathcal{G}$. Namely, it is sufficient to consider a model family reduced to a singleton and to disregard the penalty. This is also a situation where the choice of $\mathfrak{G}$ can be of interest. If we make the choice $\mathfrak{G} = \mathcal{G}$, then the target $g^*_{\mathcal{G}}$ is the best available function in the model $\mathcal{G}$. In this case, the bias term of the bound vanishes. By adding to the left and right of the obtained inequality the quantity $L(g^*_{\mathcal{G}}, g^*)$, where $g^*$ is the minimum average loss function over a larger class [e.g., $L^2(P)$], it is then possible to obtain a constant 1 in front of the bias term (instead of $\frac{K+1}{K-1} > 1$). However, this does not come completely for free since we must consider $g^*_{\mathcal{G}}$ instead of $g^*$ when checking for assumption (H3). This assumption may actually be harder to check for in practice, because usually $g^*$ has a simple, closed form (e.g., the Bayes classifier in a classification framework), whereas $g^*_{\mathcal{G}}$ depends on the approximation properties of model $\mathcal{G}$. Under certain convexity assumptions of the risk and of the model, it was shown in [4] that (H3) holds in this setting; this way we retrieve a bound in all points similar to single-model ERM results of [4].

4.2. *Application to support vector machines.* We now expose briefly the key elements needed to apply Theorem 4.3 to the SVM framework. Remember that in the case of SVMs, the natural loss function to consider is the hinge loss function $\ell(g) = (x, y) \mapsto (1 - yg(x))_+$: this is the empirical loss which is minimized (subject to regularization) to find a classifier $\widehat{g}$. Interpreting the SVM procedure as a penalized model selection procedure (see Section 2.3), we intend to apply Theorem 3.1.

To this end, we first discretize the continuous family of models $(\mathcal{B}(R))_{R \in \mathbb{R}}$ over a certain family of values of the radii: thus, our collection of models will be $(\mathcal{B}(R))_{R \in \mathcal{R}}$, where $\mathcal{R}$ is an appropriate discrete set of positive real numbers. We now have to check assumptions (H1)–(H4) of Theorem 4.3. The detailed analysis is exposed in Section 6.3 and the following statement sums up the obtained results:

THEOREM 4.4. *Let $\mathcal{R}$ be a countable set of positive real numbers, $\mathfrak{G} = L^2(P_X)$, and $\ell$ the hinge loss function.*

*In setting (S1) under assumptions (A1), (A2) and (A3), the family of models $(\mathcal{B}(R))_{R \in \mathcal{R}}$ satisfies hypotheses (H1) to (H4) of Theorem 4.3 with the following parameter values:*

$$b_R = 1 + MR; \qquad C_R = 2\left(\frac{MR}{\eta_1} + \frac{1}{\eta_0}\right);$$

$$r_R^* \le 16\frac{C_R^2}{\sqrt{n}}\inf_{d \in \mathbb{N}}\left(\frac{d}{\sqrt{n}} + \frac{\eta_1}{M}\sqrt{\sum_{j>d}\lambda_j}\right).$$

*In setting (S2) under assumptions (A1) and (A2), the family of models $(\mathcal{B}(R))_{R \in \mathcal{R}}$ satisfies hypotheses (H1) to (H4) of Theorem 4.3 with the following parameter values:*

$$b_R = 1 + MR; \qquad C_R = \left(MR + \frac{1}{\eta_0}\right); \qquad r_R^* \le 2500M^{-2}C_R^2 x_*^2(n),$$

*where $x_*$ is as in the definition of setting (S2).*

Once assumptions (H1)–(H4) are granted, the remaining task in order to prove Theorem 3.1 is to formalize precisely how to back and forth between the continuous regularization and the discrete sets of models $(\mathcal{B}(R))_{R \in \mathcal{R}}$. The details are given in Section 6.4.

## 5. Discussion and conclusion.

5.1. *Relation to other work.* In this section we compare our result to earlier work. The properties of the generalization error of the SVM algorithm have been investigated in various ways (we omit here the vast literature on algorithmic aspects of the SVM with which the present paper is not concerned). To this regard, we distinguish between two types of results: the first type are error bounds. They bound the difference between the empirical and true expected loss of an estimator. The second type are excess loss inequalities which relate the risk of the estimator to the Bayes risk.

5.1.1. *Error bounds.* The first result about the SVM algorithm is due to Vapnik; who proved that the fat-shattering dimension (see, e.g., [1] for a definition) at scale 1 of the set $\{(x, y) \mapsto y\langle k(x, \cdot), f\rangle_k + b = yf(x) + b : f \in \mathcal{H}_k, \|f\|_k \leq R, b \in \mathbb{R}\}$ on a sample $X_1, \ldots, X_n$ is bounded by $D^2 R^2$, where $D$ is the radius of the smallest ball enclosing the sample in feature space, which can be computed as $D = \inf_{g \in \mathcal{H}_k} \max_{i=1,\ldots,n} \|k(X_i, \cdot) - g\|_k$ or, equivalently, $D^2 := \max_{\substack{\|\beta\|_1 \leq 1 \\ \beta_i \geq 0}} \sum_{i=1}^{n} \beta_i k(X_i, X_i) - \sum_{i,j} \beta_i \beta_j k(X_i, X_j)$.

This bound is known as the "radius-margin" bound since it involves the ratio of the radius of the sphere enclosing the data in feature space and of the (geometrical) margin of separation of the data which is equal to $1/R$ when the scaling is chosen such that the points lying on the margin (the "support vectors") have output value in $\{-1, 1\}$.

The first formal error bounds on large margin classifiers were proven by Bartlett [2]. In these bounds, the misclassification error $\mathcal{E}(f)$ of a real-valued classifier $f$ is compared to the fraction of the sample which are misclassified or almost misclassified, that is, which have margin less than a certain (positive) value. In later work, it was noticed that for classes of functions such as $\mathcal{B}(R)$, the spectrum of the kernel operator [27] plays an important role in capacity analysis.

More recent bounds on the capacity of such classes, involving Rademacher averages, have confirmed this role. We reproduce here a particularly elegant bound based on this technique (Theorem 21 of [5], slightly adapted for our notation):

THEOREM 5.1. *Let $R > 0$; for any $x > 0$, with probability at least $1 - 4e^{-x}$, for all $f \in \mathcal{B}(R)$,*

$$P\theta(f) \leq P_n[\ell(f) \wedge 1] + \frac{4R}{\sqrt{n}} \sqrt{\frac{1}{n} \sum_{i=1}^{n} k(X_i, X_i)} + 9\sqrt{\frac{x}{2n}}.$$

Error bounds as the above are typically valid for any function in $\mathcal{B}(R)$ uniformly. They thus do not take into account the specificity of the SVM algorithm. Also, for an error bound, we cannot expect a better convergence rate than $n^{-1/2}$ of the empirical loss to the true average loss, since for a single function this is the rate given asymptotically by the central limit theorem.

The term $\sum_{i=1}^{n} k(X_i, X_i)$ in the above theorem is the trace of the so-called Gram matrix (matrix of inner products of the data points in feature space). Its expected value under the sampling of the data is precisely $n$ times the trace of the kernel operator, that is, the sum of its eigenvalues. If we compare this to our main result Theorem 3.1, in setting (S1), we see that our complexity penalty is always of smaller order (up to a constant factor, and to the relation between empirical and true spectrum, which we do not cover here, but is studied, e.g., in [8, 28]).

In a different direction, in [11] are presented error bounds for regularization algorithms which explicitly involve the regularization parameter.

5.1.2. *Excess loss inequalities.* Studying the behavior of relative (or excess) loss has been at the heart of recent work in the statistical learning field. Some results have been developed specifically for regularization algorithms of the type (2.2). In particular, asymptotic results on the consistency of the SVM algorithm, that is, convergence of the risk toward Bayes risk, were obtained by Steinwart in [31].

Using a leave-one-out analysis of the SVM algorithm and techniques similar to those in [11], Zhang [40] obtained sharp bounds on the difference between the risk of the SVM classifier and the Bayes risk of the form

$$\mathbb{E}[\ell(f_n)] - c\mathbb{E}[\ell(f^*)],$$

where $c > 1$. However, because of this last strict inequality, this means that one cannot directly obtain information about the convergence $L(f_n, f^*)$ to zero from these results as soon as $\mathbb{E}\ell(f^*)$ is nonzero.

Studying the convergence of $L(f_n, f^*)$ opens the door to complexity penalties that decrease faster than $n^{-1/2}$, because the final goal is to compare directly the true average loss of the target and the estimated function, not their empirical loss. The so-called "localized approach" (that we followed in this paper) is a theoretical device used to prove such improved rates. Introduced in the statistical community for the general study of $M$-estimation, it has become widespread recently in the learning theoretical community; see, for example, [3, 4, 6, 16, 17, 20, 23, 25].

Concerning more specifically the SVM, recent works have concentrated on obtaining faster rates of convergence in various senses. In [12], the $q$-soft margin SVM is studied (i.e., when the considered loss function is $\ell^q$ for $q > 1$). In [26], the SVM is studied from the point of view of inverse problems. In [32], convergence properties of the standard SVM is studied in the case of the Gaussian kernel. In the above references, to obtain the best bounds on the rates of convergence, the regularization parameter $\Lambda_n$ (and, in the latter reference, the width of the Gaussian kernel) must have a prescribed decrease as a function of the sample size $n$, depending on a priori knowledge on regularity properties of the function $f^*$ (or $\eta$). Therefore, these results do not display adaptivity with respect to the regularity of $f^*$.

In the recent paper [33], a general inequality for regularized risk minimizers was derived, applying, in particular, to the SVM framework. The main differences in this work with respect to our framework are the following:

- a general family of possible loss functions (which includes hinge loss and square loss) is considered;

- a general condition on the loss and the generating probability distribution is considered, covering, in particular, the general Tsybakov's noise setting for classification (but without adaptivity to this regard);
- the regularization considered is fixed to be the squared RKHS norm;
- the capacity of the kernel space is measured in terms of universal $L^2$ entropy.

While our work has obviously less generality concerning the first two points, our results are sharper concerning the two last ones. One of our main goals here was to study precisely what was the minimal order of the penalty with which we could prove an oracle inequality for the loss function used in the SVM. Furthermore, our setting (S2) relies on a capacity measure of the kernel space based on the spectral properties of the associated integral operator, which is sharper than universal entropy in this setting. Again, the approach we followed here was inspired by an analogy of the SVM with the more classical regularized least squares regression, which is by now relatively well understood, and where the optimal results concerning the two last above points are known to be sharper than those obtained in [33]. Our investigation was driven by the question of how much of these precise results could be carried over to the SVM setting.

Finally, while our results demonstrate the adaptivity of support vector machines with respect to the approximation properties by the RKHS $\mathcal{H}_k$ of the target $f*$, we do not tackle the question of full adaptivity with respect to Tsybakov's noise condition. Only recently have results been obtained in this direction [18, 34, 36].

5.2. *Conclusion.* Summing up our findings, we have brought forth a general theorem allowing to derive oracle inequalities for penalized model selection methods. Application of this theorem to support vector machines has led to precise sufficient conditions for the form of the regularization function to be used in order to obtain oracle inequalities for the hinge loss. In particular, under the assumptions considered here about the probability distribution $P(Y|X)$, the bound we obtain gets better if we use a linear regularizer in the Hilbert norm rather than the standard quadratic one.

This result thus brings forth the interesting question of whether a SVM-type algorithm using a lighter (linear in the Hilbert norm) regularizer would yield improved practical results. Several issues are in play here. First a practical issue: a disadvantage of a linear regularizer is that the associated optimization problem, although convex, is not as easily tractable from an algorithmic point of view as the squared-norm regularization. Second, a theoretical issue, namely, whether a corresponding lower bound holds, which would prove that the linear regularizer is indeed better. This is the case for regularized least squares in the Gaussian noise; for the SVM, lower bounds remain very largely an open problem. And third, a crucial issue both theoretical and practical, and not tackled here, is that the multiplicative factor $\Lambda_n$ in (2.2) is seldom taken equal to some a priori fixed function of $n$ in practice. Instead, it is typically picked by cross-validation. It is important to bring into focus the fact that, even if the quadratic regularizer was suboptimal for a fixed penalty scheme, this may still be compensated by the cross-validation step for the multiplicative factor $\Lambda_n$, which could implicitly "correct" this effect. We believe this issue has not been studied in current work on SVMs, and that it is a central point to be studied in the future in order to reconciliate theory and practice.

Several other mathematical problems remain open. Ideally, one would hope to obtain the same kind of result for the full SVM algorithm instead of the SVM$_0$ considered here. We mentioned in our comments after the main theorem a possible extension from our "gap" condition to a general Tsybakov's noise condition. This would give rise to an additional term for which we cannot always ensure that it is only a negligible remainder as the sample size grows to infinity. Therefore, the question of full adaptivity to Tsybakov's noise remains generally open. Finally, it is not clear whether our sufficient minimum rate conditions for the penalty are minimal: it would be interesting to investigate whether a lower order penalty would, for example, yield an inconsistent estimator.

## 6. Proofs.

6.1. *Proof of Lemma* 2.1.   We start with proving (i). We can write

$$\mathbb{E}[\ell(g)] = \mathbb{E}[\eta(X)(1 - g(X))_+ \\ + (1 - \eta(X))(1 + g(X))_+].$$

We will prove that, for each fixed $x$, $s^*(x)$ minimizes the expression in the expectation. Let's study the function $g \mapsto \eta(1-g)_+ + (1-\eta)(1+g)_+$. It is easy to see that for $\eta \in [\frac{1}{2}, 1]$ it is minimized for $g = 1$, and for $\eta \in [0, \frac{1}{2}]$ it is minimized for $g = -1$. This means that, in all cases, the minimum is reached at $g = s^*$. Finally, it is easy to see that this minimum is unique whenever $\eta \notin \{0, \frac{1}{2}, 1\}$, hence, $f^* = s^*$ a.s. on this set. (Notice additionally that, for $\eta = 1$, any $g \geq 1$ reaches the minimum, for $\eta = 0$, any $g \leq -1$ reaches the minimum and for $\eta = \frac{1}{2}$, any $g \in [-1, 1]$ reaches the minimum.)

We now turn to (ii). Considering (i), we can arbitrarily choose $f^* = s^*$. We then have to prove that

$$\mathbb{E}[\mathbf{1}\{Yg(X) \leq 0\} - \mathbf{1}\{Ys^*(X) \leq 0\}] \\ \leq \mathbb{E}[(1 - Yg(X))_+ - (1 - Ys^*(X))_+].$$

We know that the right-hand side is nonnegative. Moreover, the random variable in the left-hand side is positive (and thus equal to 1) if and only if $Yg(X) \leq 0$ and $Ys^*(X) \geq 0$, in which case $(1 - Yg(X))_+ \geq 1$ and $(1 - Ys^*(X))_+ = 0$ (since $s^*$ takes its values in $\{-1, 1\}$). This proves the inequality.

6.2. *Proof of Theorem* 4.3. To prove Theorem 4.3, we first state the key technical result concerning a localized uniform control of an empirical process.

THEOREM 6.1. *Let $\mathcal{F}$ be a class of measurable, square integrable functions such that for all $f \in \mathcal{F}$, $Pf - f \leq b$. Let $w(f)$ be a nonnegative function such that $\mathsf{Var}[f] \leq w(f)$. Let $\phi$ be a sub-root function, $D$ be some positive constant and $r^*$ be the unique positive solution of $\phi(r) = r/D$. Assume that the following holds:*

$$(6.1) \qquad \forall r \geq r^* \qquad \mathbb{E}\left[0 \vee \left(\sup_{f \in \mathcal{F}: w(f) \leq r}(P - P_n)f\right)\right] \leq \phi(r).$$

*Then, for all $x > 0$ and all $K > D/7$, the following inequality holds with probability at least $1 - e^{-x}$:*

$$\forall f \in \mathcal{F} \qquad Pf - P_n f \leq K^{-1}w(f) + \frac{50K}{D^2}r^* + \frac{(K + 9b)x}{n}.$$

*If additionally, the convex hull of $\mathcal{F}$ contains the null function, the same is true when the positive part in (6.1) is removed.*

Note that this result is very similar to Theorem 3.3 in [3] which was obtained using techniques from [21]. We use similar techniques to obtain the version presented here.

We will need to transform assumption (6.1), using the following technical lemma which is a form of the so-called "peeling device"; the version presented here is very close to a similar lemma in [22].

LEMMA 6.2. *If $\phi$ is a sub-root function such that for any $r \geq r^* \geq 0$,*

$$\mathbb{E}\left[0 \vee \left(\sup_{f \in \mathcal{F}: w(f) \leq r} Pf - P_n f\right)\right] \leq \phi(r),$$

*one has for any $r \geq r^*$,*

$$\mathbb{E}\left[\sup_{f \in \mathcal{F}} \frac{Pf - P_n f}{w(f) + r}\right] \leq 4\frac{\phi(r)}{r},$$

*and when $0 \in \mathrm{conv}\mathcal{F}$, the same is true if the positive part is removed in the previous condition.*

PROOF. We choose some $x > 1$. In the calculations below a supremum over an empty set is considered as 0. We have

$$\sup_{f \in \mathcal{F}} \frac{Pf - P_n f}{w(f) + r}$$

$$\leq \sup_{f \in \mathcal{F}: w(f) \leq r} \frac{(Pf - P_n f)_+}{w(f) + r} + \sum_{k \geq 0} \sup_{f \in \mathcal{F}: rx^k \leq w(f) \leq rx^{k+1}} \frac{(Pf - P_n f)_+}{w(f) + r}$$

$$\leq \frac{1}{r} \sup_{f \in \mathcal{F}: w(f) \leq r} (Pf - P_n f)_+ + \sum_{k \geq 0} \sup_{f \in \mathcal{F}: rx^k \leq w(f) \leq rx^{k+1}} \frac{(Pf - P_n f)_+}{rx^k + r}$$

$$\leq \frac{1}{r} \left( \sup_{f \in \mathcal{F}: w(f) \leq r} (Pf - P_n f)_+ + \sum_{k \geq 0} \sup_{f \in \mathcal{F}: w(f) \leq rx^{k+1}} \frac{(Pf - P_n f)_+}{1 + x^k} \right).$$

In the general case, note that $\sup_{a \in A}(0 \vee a) = 0 \vee \sup_{a \in A} a$. In the case where $\operatorname{conv} \mathcal{F}$ contains the null function, one has $\sup_{f \in \mathcal{F}} Pf - P_n f = \sup_{f \in \operatorname{conv} \mathcal{F}} Pf - P_n f \geq 0$ so that $\sup_{f \in \mathcal{F}}(Pf - P_n f)_+ = \sup_{f \in \mathcal{F}} Pf - P_n f$, which allows us to remove the positive part in the condition for $\phi$.

So, taking the expectation, we obtain

$$\mathbb{E}\left[\sup_{f \in \mathcal{F}} \frac{Pf - P_n f}{w(f) + r}\right] \leq \frac{1}{r}\left( \phi(r) + \sum_{k \geq 0} \frac{\phi(rx^{(k+1)})}{1 + x^k} \right)$$

$$\leq \frac{\phi(r)}{r}\left( 1 + \sum_{k \geq 0} \frac{x^{(k+1)/2}}{1 + x^k} \right)$$

$$\leq \frac{\phi(r)}{r}\left( 1 + x^{1/2}\left( \frac{1}{2} + \sum_{k \geq 1} x^{-k/2} \right) \right)$$

$$\leq \frac{\phi(r)}{r}\left( 1 + x^{1/2}\left( \frac{1}{2} + \frac{1}{x^{1/2} - 1} \right) \right),$$

where we have used the sub-root property for the second inequality. It is then easy to check that the minimum of the right-hand side is attained at

$$x = (1 + \sqrt{2})^2.$$

Plugging this value in the right-hand side, we obtain the result. □

PROOF OF THEOREM 6.1. The main technical tool of the proof is Talagrand's concentration inequality (here we use an improved version proved in [10]). We recall it briefly as follows.

Let $X_i$ be independent variables distributed according to $P$, and $\mathcal{F}$ a set of functions from $\mathcal{X}$ to $\mathbb{R}$ such that $\mathbb{E}[f] = 0$, $\|f\|_\infty \leq c$ and $\operatorname{Var}[f] \leq \sigma^2$ for any $f \in \mathcal{F}$. Let $Z = \sup_{f \in \mathcal{F}} \sum_{i=1}^{n} f(X_i)$. Then with probability $1 - e^{-x}$, it holds that

$$(6.2) \qquad Z \leq \mathbb{E}Z + \sqrt{2x(n\sigma^2 + 2c\mathbb{E}[Z])} + \frac{cx}{3}.$$

We will apply this inequality to the rescaled set of functions

$$\mathcal{F}_r = \left\{ \frac{Pf - f}{w(f) + r}, f \in \mathcal{F} \right\},$$

where we assume $r \geq r^*$. The precise choice for $r$ will be decided later. We now check the assumptions on the supremum norm and the variance of functions in $\mathcal{F}_r$. We have

$$\sup_{f \in \mathcal{F}} \sup_{x \in \mathcal{X}} \frac{Pf - f(x)}{r + w(f)} \leq \frac{b}{r};$$

and, recalling the hypothesis that $\mathsf{Var}[f] \leq w(f)$, the following holds:

$$\mathsf{Var}\left[ \frac{f(X)}{w(f) + r} \right] \leq \frac{w(f)}{(w(f) + r)^2} \leq \frac{w(f)}{4rw(f)} = r^{-1}/4,$$

where we have used the fact that $2ab \leq a^2 + b^2$. Introducing the following random variable

$$(6.3) \qquad V_r = \sup_{f \in \mathcal{F}} \frac{Pf - P_n f}{w(f) + r},$$

we thus obtain by application of (6.2) that, with probability at least $1 - e^{-x}$,

$$(6.4) \qquad V_r \leq \mathbb{E}[V_r] + \sqrt{\frac{x}{2rn} + \frac{4xb\mathbb{E}[V_r]}{rn}} + \frac{xb}{3rn}.$$

It follows from Lemma 6.2 that $\mathbb{E}[V_r] \leq 4\phi(r)/r$. Plugging this into (6.4), and recalling that $r^*$ is the unique solution of $\phi(r) = r/D$, we obtain that, for all $x > 0$, and $r \geq r^*$, the following inequality hold with probability at least $1 - e^{-x}$:

$$(6.5) \qquad \begin{aligned} &\forall f \in \mathcal{F} \\ &\frac{Pf - P_n f}{w(f) + r} \leq \inf_{\alpha > 0} \left( 4\frac{1+\alpha}{D}\sqrt{\frac{r^*}{r}} + \sqrt{\frac{x}{2nr}} + \left( \frac{1}{3} + \frac{1}{\alpha} \right)\frac{bx}{rn} \right). \end{aligned}$$

Here, we have used the fact that, for $r \geq r^*$, $\phi(r)/r \leq \sqrt{r^*/rD^2}$ and that $2\sqrt{ab} \leq \alpha a + b/\alpha$.

Now given some constant $K$, we want to find $r \geq r^*$ such that $V_r \leq 1/K$ (with high probability). This corresponds to finding $r$ such that the left-hand side of (6.5) is upper bounded by $1/K$.

Denote $A_1 = 4(1 + \alpha)\sqrt{r^*}/D + \sqrt{x/2n}$ and $A_2 = (1/3 + 1/\alpha)bx/n$. Then we have to find $r$ such that $A_1 r^{-1/2} + A_2 r^{-1} \leq K^{-1}$. It can be easily checked that this is satisfied if

$$(6.6) \qquad\qquad r \geq K^2 A_1^2 + 2A_2 K.$$

We have

$$K^2 A_1^2 + 2A_2 K \leq 32(1 + \alpha)^2 \frac{K^2 r^*}{D} + \frac{x}{n}(K^2 + 2bK/3 + 2bK/\alpha).$$

Taking $\alpha = 1/4$, we conclude that (6.6) is satisfied when the following holds:

$$r \geq 50 \frac{K^2}{D^2} r^* + (K^2 + 9bK)\frac{x}{n}.$$

Note that $K > D/7$ ensures that the lower bound above is greater than $r^*$. We can thus take $r$ equal to this value.

Combining the above results concludes the proof of Theorem 6.1. $\square$

We are now in a position to proceed to the proof of the main model selection theorem.

PROOF OF THEOREM 4.3. The main use of hypotheses (H1), (H2) and (H4) will be to apply Theorem 6.1 to the class

$$\mathcal{F}_{m,g_0} = \{\ell(g) - \ell(g_0), g \in \mathcal{G}_m\}$$

for some $m \in \mathcal{M}, g_0 \in \mathcal{G}_m$ with the choice $w(f) = \min\{d^2(g, g_0) | g \in \mathcal{G}_m, \ell(g) - \ell(g_0) = f\}$, so that, using hypotheses (H1), (H2), (H4) and the fact that the null function belongs to the class, we obtain that, for any arbitrary $K > C/7$, with probability at least $1 - e^{-x}$,

$$(6.7) \qquad \forall g \in \mathcal{G}_m$$
$$(P - P_n)(\ell(g) - \ell(g_0)) \leq K^{-1}d^2(g, g_0) + \frac{50K}{C^2}r^* + \frac{(K + 9b)x}{n}.$$

For each $m \in \mathcal{M}$, we define $u_m$ and $g_m$ as functions in $\mathcal{G}_m$ satisfying, respectively,

$$\begin{cases} d(u_m, g^*) = \inf_{g \in \mathcal{G}_m} d(g, g^*), \\ L(g_m, g^*) = \inf_{g \in \mathcal{G}_m} L(g, g^*). \end{cases}$$

[If these infima are not attained, one can choose $u_m, g_m$ such that $d(u_m, g^*)$, $L(g_m, g^*)$ are arbitrary close to the inf, and use a dominated convergence argument at the end of the proof.]

Now, for any $m \in \mathcal{M}, g_m \in \mathcal{F}_m$,

$$
\begin{aligned}
(6.8) \quad & L(\widetilde{g}, g^*) - L(g_m, g^*) \\
&= P\ell(\widetilde{g}) - P\ell(g_m) \\
&= P_n\ell(\widetilde{g}) - P_n\ell(g_m) + (P - P_n)(\ell(\widetilde{g}) - \ell(g_m)) \\
&\leq \operatorname{pen}(m) - \operatorname{pen}(\widetilde{m}) + \rho_m + (P - P_n)(\ell(\widetilde{g}) - \ell(g_m)),
\end{aligned}
$$

where the last inequality stems from the definition of $\widetilde{g}$.

Denoting $\widetilde{m}$ the model containing $\widetilde{g}$, we decompose the last term above:

$$
\begin{aligned}
(6.9) \quad (P - P_n)(\ell(\widetilde{g}) - \ell(g_m)) &= (P - P_n)(\ell(\widetilde{g}) - \ell(u_{\widetilde{m}})) \\
&\quad + (P - P_n)(\ell(u_{\widetilde{m}}) - \ell(g_m)).
\end{aligned}
$$

We will bound both terms separately. For the first term, we use (6.7): for any $m' \in \mathcal{M}$ and an arbitrary $K_{m'} > C_{m'}/7$, with probability at least $1 - e^{-x_{m'} - \xi}$, for all $g \in \mathcal{G}_{m'}$ we have

$$
\begin{aligned}
(6.10) \quad & (P - P_n)(\ell(g) - \ell(u_{m'})) \\
&\leq K_{m'}^{-1} d^2(g, u_{m'}) + \frac{50 K_{m'}}{C_{m'}^2} r_{m'}^* \\
&\quad + \frac{(K_{m'} + 9b_{m'})(x_{m'} + \xi)}{n}.
\end{aligned}
$$

By the union bound, this inequality is valid simultaneously for all $m' \in \mathcal{M}$ with probability $1 - e^{-\xi}$, so that it holds, in particular, for $m' = \widetilde{m}, g = \widetilde{g}$ with this probability. Finally, note that, for $g \in \mathcal{G}_{\widetilde{m}}$,

$$
(6.11) \quad d^2(g, u_{\widetilde{m}}) \leq (d(g, g^*) + d(u_{\widetilde{m}}, g^*))^2 \leq 4d^2(g, g^*).
$$

For the second term of (6.9), we will use the following Bernstein inequality: for any $m_1, m_2 \in \mathcal{M}$, we have, with probability $1 - \exp(-x_{m_1} - x_{m_2} - \xi)$,

$$
\begin{aligned}
(6.12) \quad & (P - P_n)(\ell(u_{m_1}) - \ell(g_{m_2})) \\
&\leq \sqrt{2(x_{m_1} + x_{m_2} + \xi) \frac{\mathsf{Var}[\ell(u_{m_1}) - \ell(g_{m_2})]}{n}} \\
&\quad + \frac{\max(b_{m_1}, b_{m_2})(x_{m_1} + x_{m_2} + \xi)}{6n}.
\end{aligned}
$$

Now, using assumption (4.1), if $b_{m^*} = \max(b_{m_1}, b_{m_2})$,

$$
\max(b_{m_1}, b_{m_2})(x_{m_1} + x_{m_2}) \leq 2b_{m^*} x_{m^*} \leq 2b_{m_1} x_{m_1} + 2b_{m_2} x_{m_2}.
$$

We now deal with the first term of the bound (6.12): for any $g \in \mathcal{G}_{m_1}$,

$$\sqrt{2(x_{m_1} + x_{m_2} + \xi) \frac{\mathsf{Var}[\ell(u_{m_1}) - \ell(g_{m_2})]}{n}}$$

$$\leq \sqrt{4(x_{m_1} + x_{m_2} + \xi) \frac{(d^2(u_{m_1}, g^*) + d^2(g_{m_2}, g^*))}{n}}$$

$$\leq 2\sqrt{(x_{m_1} + x_{m_2} + \xi) \frac{d^2(g, g^*)}{n}} + 2\sqrt{(x_{m_1} + x_{m_2} + \xi) \frac{d^2(g_{m_2}, g^*)}{n}}$$

$$\leq K_{m_1}^{-1} d^2(g, g^*) + K_{m_2}^{-1} d^2(g_{m_2}, g^*) + \frac{(K_{m_1} + K_{m_2})(x_{m_1} + x_{m_2} + \xi)}{n},$$

where the first inequality follows from hypothesis (H2) followed by the triangle inequality, and the second from the definition of $u_{m_1}$. Anticipating somewhat the end of the proof, we will choose $K_m = \alpha C_m$ for some fixed $\alpha$, so that, using again assumption (4.1) like above, it is true that

$$(K_{m_1} + K_{m_2})(x_{m_1} + x_{m_2}) \leq 4K_{m_1} x_{m_1} + 4K_{m_2} x_{m_2}.$$

Therefore, (6.12) becomes, with probability $1 - \exp(-x_{m_1} - x_{m_2} - \xi)$, for any $g \in \mathcal{G}_{m_1}$,

$$\begin{aligned}(6.13) \quad & (P - P_n)(\ell(u_{m_1}) - \ell(g_{m_2})) \\ & \leq K_{m_1}^{-1} d^2(g, g^*) + K_{m_2}^{-1} d^2(g_{m_2}, g^*) \\ & \quad + \frac{(12K_{m_1} + b_{m_1})(x_{m_1} + \xi)}{3n} \\ & \quad + \frac{(12K_{m_2} + b_{m_2})(x_{m_2} + \xi)}{3n}. \end{aligned}$$

Bound, (6.13) is therefore valid for all $m_1, m_2 \in \mathcal{M}$ simultaneously with probability $1 - \exp(-\xi)$, and, in particular, for $m_1 = \widetilde{m}, m_2 = m, g = \widetilde{g}$.

Putting together (6.9), (6.10), (6.11) and (6.13), we obtain that, with probability $1 - 2\exp(-\xi)$, for all $m \in \mathcal{M}$,

$$\begin{aligned}(6.14) \quad & (P - P_n)(\ell(\widetilde{g}) - \ell(g_m)) \\ & \leq 5K_{\widetilde{m}}^{-1} d^2(\widetilde{g}, g^*) + K_m^{-1} d^2(g_m, g^*) + \frac{50K_{\widetilde{m}}}{C_{\widetilde{m}}^2} r_{\widetilde{m}}^* \\ & \quad + \frac{(15K_{\widetilde{m}} + 28b_{\widetilde{m}})(x_{\widetilde{m}} + \xi)}{3n} + \frac{(12K_m + b_m)(x_m + \xi)}{3n}. \end{aligned}$$

Now choosing $K_m = 5KC_m$ (note that we have $K_m > C_m/7$ as required, since $K > 1$), and replacing $\xi$ by $\xi + \log(2)$, recalling inequality (6.8) and the hypothesis (4.2) on the penalty function, we thus obtain that, with probability $1 - \exp(-\xi)$, for any $m \in \mathcal{M}$,

$$
\begin{aligned}
L(\widetilde{g}, g^*) &- L(g_m, g^*) \\
&\leq \operatorname{pen}(m) - \operatorname{pen}(\widetilde{m}) + C_{\widetilde{m}}^{-1} K^{-1} d^2(\widetilde{g}, g^*) + \tfrac{1}{5} C_m^{-1} K^{-1} d^2(g_m, g^*) \\
&\quad + \operatorname{pen}(\widetilde{m}) + \operatorname{pen}(m) + \rho_m \\
&\leq K^{-1} L(\widetilde{g}, g^*) + \tfrac{1}{5} K^{-1} L(g_m, g^*) + 2\operatorname{pen}(m) + \rho_m,
\end{aligned}
$$

using hypothesis (H4). This leads to the conclusion for the deviation inequality of the model selection theorem.

For the inequality in expected risk, we go back to inequality (6.14), with the choice $K_m = 5KC_m$; also using (6.8), we conclude that, for any $\xi > 0$, the following inequality holds with probability $1 - \exp(-\xi)$:

$$
\begin{aligned}
(6.15) \qquad L(\widetilde{g}, g^*) &- L(g_m, g^*) \\
&\leq K^{-1} L(\widetilde{g}, g^*) + \tfrac{1}{5} K^{-1} L(g_m, g^*) + \operatorname{pen}(m) - \operatorname{pen}(\widetilde{m}) + \rho_m \\
&\quad + \frac{250 K C_{\widetilde{m}}}{D_{\widetilde{m}}^2} r_{\widetilde{m}}^* + \frac{B_{\widetilde{m}}(x_{\widetilde{m}} + \xi + \log(2))}{3n} \\
&\quad + \frac{B_m(x_m + \xi + \log(2))}{3n}.
\end{aligned}
$$

The point is now to linearize the products $B_m \xi$. To do so, we use the following Young's inequality valid for any positive $x, y$:

$$
xy \leq \exp\left(\frac{x}{2}\right) + 2y \log y,
$$

with $x = \xi, y = B_m$, so that, putting $u = \exp(\xi/2)$, and using the hypothesis (4.4) on the penalty function, we obtain that, with probability $1 - (u^{-2} \wedge 1)$,

$$
\begin{aligned}
(6.16) \qquad L(\widetilde{g}, g^*) - L(g_m, g^*) &\leq K^{-1} L(\widetilde{g}, g^*) + \tfrac{1}{5} K^{-1} L(g_m, g^*) \\
&\quad + 2\operatorname{pen}(m) + \rho_m + \frac{u}{n}.
\end{aligned}
$$

Integrating concludes the proof. $\square$

6.3. *Proof of Theorem* 4.4. The purpose of Theorem 4.4 is to check that conditions (H1) to (H4) of the general model selection Theorem 4.3 are satisfied for settings (S1) and (S2) of the SVM. We will split the proofs into several results corresponding to the different hypotheses.

LEMMA 6.3. *Under assumption* (A1), *hypothesis* (H1) *is satisfied for* $b_R = MR + 1$.

Proof. We use the reproducing property of the kernel to conclude that

$$\forall g \in \mathcal{B}(R) \qquad |\ell(yg(x))| \leq 1 + |g(x)|$$
$$= 1 + |\langle g, k(x, \cdot) \rangle|_k$$
$$\leq 1 + \|g\|_k \|k(x, \cdot)\|_k$$
$$= 1 + \|g\|_k \sqrt{k(x,x)} \leq 1 + MR. \qquad \square$$

We now check conditions (H2) and (H3). This differs according to the setting, because we make a different choice for the pseudo-distance $d$ depending on the setting considered.

Lemma 6.4 [Setting (S1)]. *Under assumptions* (A1), (A2) *and* (A3), *conditions* (H2) *and* (H3) *of Theorem* 4.3 *are satisfied for the choice*

$$d_1(g, g') = \mathbb{E}[(g - g')^2]; \qquad C_R = 2\left(\frac{MR}{\eta_1} + \frac{1}{\eta_0}\right).$$

Proof. Obviously, (H2) is satisfied since $\ell$ is a Lipschitz function, so that $|\ell(yg(x)) - \ell(yg'(x))| \leq |g(x) - g'(x)|$.

We will obtain the result we look for if we can bound uniformly in $x$ the ratio

$$\frac{\mathbb{E}[(g - s^*)^2 \mid X = x]}{\mathbb{E}[\ell(g) - \ell(s^*) \mid X = x]}.$$

Remember that, for $g \in \mathcal{B}(R)$, the reproducing property of the kernel and assumption (A1) imply $\|g\|_\infty \leq M\|g\|_k \leq MR$. Let us consider without loss of generality the case $s^* = 1$ (i.e., $\eta = \mathbb{P}[Y = 1 | X = x] \geq \frac{1}{2}$). We then have to bound the ratio

$$\frac{(1 - g)^2}{\eta(1 - g)_+ + (1 - \eta)(1 + g)_+ - 2(1 - \eta)}.$$

For $g \leq -1$, this becomes $\frac{(1-g)^2}{\eta(1-g) - 2(1-\eta)}$; putting $x = -g - 1 \in [0, MR]$, this can be rewritten as

$$\frac{(x + 2)^2}{\eta x + 2(2\eta - 1)} \leq \frac{2x^2 + 8}{\eta x + 2(2\eta - 1)} \leq 4MR + \frac{2}{\eta_0} \leq 2\left(\frac{MR}{\eta_1} + \frac{1}{\eta_0}\right),$$

where we have used the fact that $\eta \geq \frac{1}{2} \geq \eta_1$. For $g \geq 1$, this becomes $\frac{g-1}{1-\eta}$, which is smaller than $(MR - 1)/\eta_1$ for $g \in [1, MR]$. For $g \in [-1, 1]$, the ratio becomes $\frac{1-g}{2\eta-1}$, which is smaller than $1/\eta_0$. $\square$

LEMMA 6.5 [Setting (S2)]. *Under assumptions* (A1) *and* (A2), *conditions* (H2) *and* (H3) *of Theorem* 4.3 *are satisfied for the choice*

$$d_2(g, g') = \mathbb{E}[(\ell(g) - \ell(g'))^2]; \qquad C_R = \left(MR + \frac{1}{\eta_0}\right).$$

PROOF. Obviously, (H2) is satisfied as before. We will obtain the result we look for if we can bound uniformly in $x$ the ratio

$$\frac{\mathbb{E}[\ell(g)^2 - 2\ell(g)\ell(s^*) + \ell^2(s^*) \mid X = x]}{\mathbb{E}[\ell(g) - \ell(s^*) \mid X = x]}.$$

Notice first that

$$\mathbb{E}[\ell^2(s^*) \mid X = x] = 2\mathbb{E}[\ell(s^*) \mid X = x] = 4\min(\eta(x), 1 - \eta(x)).$$

Let us first consider the case $s^* = 1$ (i.e., $\eta \geq \frac{1}{2}$). The above ratio can be written

$$\frac{\eta(1 - g)_+^2 + (1 - \eta)(1 + g)_+((1 + g)_+ - 4) + 4(1 - \eta)}{\eta(1 - g)_+ + (1 - \eta)(1 + g)_+ - 2(1 - \eta)}.$$

For $g \leq -1$, this becomes $\frac{\eta(1 - g)^2 + 4(1 - \eta)}{\eta(1 - g) - 2(1 - \eta)}$; putting $x = -g - 1 \in [0, MR]$, this can be written as

$$\frac{\eta x^2 + 4\eta x + 4}{\eta x + 2(2\eta - 1)} = x + \frac{4 + 2x}{\eta x + 2(2\eta - 1)} \leq MR + \frac{1}{\eta_0}.$$

For $g \geq 1$, this becomes $\frac{(1 - g)^2}{g - 1} = g - 1$, which is smaller than $MR - 1$ for $g \in [1, MR]$. For $g \in [-1, 1]$, this becomes $\frac{1 - g}{2\eta - 1}$, which is smaller than $1/\eta_0$. The case $\eta < \frac{1}{2}$ can be treated in a similar way. □

Finally, we check for hypothesis (H4); this condition characterizes the complexity of the models and constitutes the meaty part of Theorem 4.4. We start with the following result which deals with setting (S1). Here we can see the (technical) reason why assumption (A3) was introduced in this setting: to relate the penalty to the spectrum of the integral operator, we use the $L^2$ distance $d_1$ as an intermediate pseudo-distance; but this requires in turn, assumption (A3) to check hypothesis (H3) (see Lemma 6.4 above).

THEOREM 6.6. *Let* $\mathcal{G}$ *be a RKHS with reproducing kernel* $k$ *such that the associated integral operator* $L_k$ *has eigenvalues* $(\lambda_i)$ *(in nonincreasing order). Let* $\ell$ *be the hinge loss function and denote* $d_1^2(g, u) = P(g - u)^2$. *Then, for all* $r > 0$ *and* $u \in \mathcal{B}(R)$,

$$\mathbb{E}\left[\sup_{\substack{g \in \mathcal{B}(R) \\ d_1^2(g, u) \leq r}} |(P - P_n)(\ell(g) - \ell(u))|\right] \leq \frac{4}{\sqrt{n}} \inf_{d \in \mathbb{N}} \left(\sqrt{dr} + 2R\sqrt{\sum_{j > d} \lambda_j}\right) := \phi_R(r).$$

*The above $\phi_R$ is a sub-root function, and the unique solution of $\phi_R(r) = r/C_R$, with $C_R \geq \eta_1^{-1}MR$, is upper bounded by*

$$r_R^* \leq 16\frac{C_R^2}{\sqrt{n}}\inf_{d\in\mathbb{N}}\left(\frac{d}{\sqrt{n}} + \frac{\eta_1}{M}\sqrt{\sum_{j>d}\lambda_j}\right).$$

To prove Theorem 6.6, we will use two technical results; the first will allow to bound the quantity we are interested in by a localized Rademacher complexity term; the second will give an upper bound on this term using the assumptions.

We introduce the following notation for Rademacher averages: let $\sigma_1, \ldots, \sigma_n$ be $n$ i.i.d. *Rademacher* random variables (i.e., such that $\mathbb{P}[\sigma_i = 1] = \mathbb{P}[\sigma_i = -1] = \frac{1}{2}$), independent of $(X_i, Y_i)_{i=1}^n$; then we define for any measurable real-valued function $f$ on $\mathcal{X} \times \mathcal{Y}$

$$(6.17) \qquad R_n f := n^{-1}\sum_{i=1}^n \sigma_i f(X_i, Y_i).$$

We then extend this notation to sets $\mathcal{F}$ of functions from $\mathcal{X} \times \mathcal{Y}$ to $\mathbb{R}$, denoting

$$R_n\mathcal{F} = \sup_{f\in\mathcal{F}} R_n f.$$

We then have the following lemma:

LEMMA 6.7. *Let $\mathcal{F}$ be a set of real functions; let $\phi$ be a 1-Lipschitz function on $\mathbb{R}$. Then for $g_0 \in \mathcal{F}$,*

$$\mathbb{E}\left[\sup_{g\in\mathcal{F}}|(P - P_n)(\phi \circ g - \phi \circ g_0)|\right] \leq 4\mathbb{E}R_n\{g - g_0 : g \in \mathcal{F}\}.$$

PROOF. By a symmetrization argument, we have

$$\mathbb{E}\left[\sup_{g\in\mathcal{F}}|(P - P_n)\phi \circ g - (P - P_n)\phi \circ g_0|\right] \leq 2\mathbb{E}\left[\sup_{g\in\mathcal{F}}|R_n(\phi \circ g - \phi \circ g_0)|\right],$$

and by symmetry of the Rademacher random variables, we have

$$\mathbb{E}\left[\sup_{g\in\mathcal{F}}|R_n(\phi \circ g - \phi \circ g_0)|\right] \leq 2\mathbb{E}\left[\sup_{g\in\mathcal{F}}(R_n(\phi \circ g - \phi \circ g_0))_+\right].$$

Since $g_0 \in \mathcal{F}$, choosing $g = g_0$, one notices that

$$\mathbb{E}\left[\sup_{g\in\mathcal{F}}(R_n(\phi \circ g - \phi \circ g_0))_+\right] = \mathbb{E}\left[\sup_{g\in\mathcal{F}}(R_n(\phi \circ g - \phi \circ g_0))\right],$$

and since $g_0$ is fixed, and $\mathbb{E} R_n \phi \circ g_0 = 0$, we obtain

$$\mathbb{E}\left[\sup_{g \in \mathcal{F}} |(P - P_n)\phi \circ g - (P - P_n)\phi \circ g_0|\right] \le 4\mathbb{E}\left[\sup_{g \in \mathcal{F}} R_n(\phi \circ g)\right].$$

Since $\phi$ is 1-Lipschitz, we can finally apply the contraction principle for Rademacher averages; then using $\mathbb{E} R_n g_0 = 0$, we obtain the result.   □

   The next lemma gives a result similar to [25], but we provide a slightly different proof (also, we are not concerned about lower bounds here). The principle of the proof below can be traced back to the work of R. M. Dudley.

Lemma 6.8.

$$\mathbb{E} R_n \{g \in \mathcal{H}_k \colon \|g\|_k \le R, \|g\|_{2,P}^2 \le r\} \le \frac{1}{\sqrt{n}} \inf_{d \in \mathbb{N}} \left(\sqrt{dr} + R\sqrt{\sum_{j > d} \lambda_j}\right)$$

$$\le \sqrt{\frac{2}{n}} \left(\sum_{j \ge 1} \min(r, R^2 \lambda_j)\right)^{1/2}.$$

Proof.   For $g \in \mathcal{H}_k$, by Lemma A.1 in the Appendix, we can decompose $g$ as

$$g(x) = \sum_{i > 0} \alpha_i \psi_i(x),$$

with $\|g\|_{2,P}^2 = \sum_{i > 0} \lambda_i \alpha_i^2$ and $\|g\|_k^2 = \sum_{i > 0} \alpha_i^2$. The above series representation holds as an equality in $\mathcal{H}_k$, and hence pointwise since the evaluation functionals are continuous in a RKHS. Let us denote

$$\Gamma(R, r) = \left\{\alpha \in \ell_2 \colon \|\alpha\|^2 \le R^2, \sum_{i > 0} \lambda_i \alpha_i^2 \le r\right\}.$$

Thus, the quantity we try to upper bound is equal to

$$\frac{1}{n}\mathbb{E}\left[\sup_{\alpha \in \Gamma(R, r)} \left|\sum_{i = 1}^{n} \sigma_i g_\alpha(X_i)\right|\right],$$

where

$$g_\alpha(X_i) = \sum_{j > 0} \alpha_j \psi_j(X_i).$$

We now write for any nonnegative integer $d$ and $\alpha \in \Gamma(R, r)$

$$\left|\sum_{i = 1}^{n} \sigma_i g_\alpha(X_i)\right| = \left|\sum_{j > 0} \alpha_j \sum_{i = 1}^{n} \sigma_i \psi_j(X_i)\right|$$

(6.18)
$$\leq \left| \sum_{j \leq d} \alpha_j \sum_{i=1}^{n} \sigma_i \psi_j(X_i) \right| + \left| \sum_{j > d} \alpha_j \sum_{i=1}^{n} \sigma_i \psi_j(X_i) \right|.$$

Applying the Cauchy–Schwarz inequality to the second term, we have

$$\left| \sum_{j > d} \alpha_j \sum_{i=1}^{n} \sigma_i \psi_j(X_i) \right| \leq \left( \sum_{j > d} \alpha_j^2 \right)^{1/2} \left( \sum_{j > d} \left( \sum_{i=1}^{n} \sigma_i \psi_j(X_i) \right)^2 \right)^{1/2}$$

$$\leq R \left( \sum_{j > d} \left( \sum_{i=1}^{n} \sigma_i \psi_j(X_i) \right)^2 \right)^{1/2}.$$

We now take the expectation with respect to $(\sigma_i)$ and $(X_i)$ in succession. We use the fact that the $(\sigma_i)$ are zero mean, uncorrelated, unity variance variables; then the fact that $\mathbb{E}[X^{1/2}] \leq \mathbb{E}[X]^{1/2}$, to obtain

$$\mathbb{E}_X \mathbb{E}_\sigma \left[ \left| \sum_{j > d} \alpha_j \sum_{i=1}^{n} \sigma_i \psi_j(X_i) \right| \right] \leq R \mathbb{E}_X \left[ \left( \sum_{j > d} \sum_{i=1}^{n} \psi_j^2(X_i) \right)^{1/2} \right]$$

$$\leq \sqrt{n} R \left( \sum_{j > d} \lambda_j \right)^{1/2},$$

where we have used the fact that $\mathbb{E}_X[\psi_j^2(X)] = \lambda_j$. We now apply exactly the same treatment to the first term of (6.18), except that we use weights $(\lambda_i)$ in the Cauchy–Schwarz inequality, yielding

$$\mathbb{E}_X \mathbb{E}_\sigma \left[ \left| \sum_{j \leq d} \alpha_j \sum_{i=1}^{n} \sigma_i \psi_j(X_i) \right| \right] \leq \left( \sum_{j \leq d} \lambda_j \alpha_j^2 \right)^{1/2} \mathbb{E}_X \left[ \sum_{j \leq d} \sum_{i=1}^{n} \lambda_j^{-1} \psi_j^2(X_i) \right]^{1/2}$$

$$\leq \sqrt{nrd}.$$

This gives the first result. The second one follows from choosing $d$ such that, for all $j > d$, $R^2 \lambda_j \leq r$, and using the inequality $\sqrt{A} + \sqrt{B} \leq \sqrt{2}\sqrt{A + B}$. $\square$

PROOF OF THEOREM 6.6. For $g$ a function $\mathcal{X} \to \mathbb{R}$, let us briefly introduce the notation $\overline{g} \colon (x, y) \in \mathcal{X} \times \mathcal{Y} \mapsto yg(x) \in \mathbb{R}$. Let us apply Lemma 6.7 to $\overline{\mathcal{F}}_u = \{\overline{g}, g \in \mathcal{F}_u\}$, where $\mathcal{F}_u = \{g \in \mathcal{H}_k : \|g\|_k \leq R, d^2(g, u) \leq r\}$. The hinge loss function $\ell$ is 1-Lipschitz, and $u \in \mathcal{F}$, hence,

$$\mathbb{E} \left[ \sup_{g \in \mathcal{F}} |(P - P_n)(\ell(g) - \ell(u))| \right] \leq 4 \mathbb{E} R_n \{\overline{g} - \overline{u}, g \in \mathcal{F}_u\} = 4 \mathbb{E} R_n \{g - u, g \in \mathcal{F}_u\},$$

where the last equality is true because of the symmetry of the Rademacher variables. Notice that, since $\|u\|_k \leq R$, we have

$$\{g - u, g \in \mathcal{F}_u\} \subset \{g - u, \|g - u\|_k \leq 2R, d^2(g, u) \leq r\};$$

since $d^2(g, u) = \mathbb{E}[(g - u)^2]$ is a norm-induced distance, we can replace $g - u$ by $g$ (by linearity) so that the above term can be upper bounded by

$$4\mathbb{E}R_n\{g \in \mathcal{H}_k : \|g\|_k \leq 2R, \|g\|_{2,P} \leq \sqrt{r}\}.$$

Using Lemma 6.8, this can be further upper bounded by

$$\frac{4}{\sqrt{n}} \inf_{d \in \mathbb{N}} \left( \sqrt{dr} + 2R \sqrt{\sum_{j > d} \lambda_j} \right),$$

which concludes the proof of the first part of the theorem.

Observe that the minimum of two sub-root functions is a sub-root function, so that $\phi_R$ is a sub-root function. We now compute an upper bound on the solution of the equation $\phi_R(r) = r/C$, which can be written

$$r = \frac{4C}{\sqrt{n}} \inf_{d \in \mathbb{N}} \left( \sqrt{rd} + 2R \sqrt{\sum_{j > d} \lambda_j} \right).$$

Notice that the infimum is a minimum since the series $\sum_{j \geq 1} \lambda_j$ is converging and, thus, the value of the right-hand side is bounded for all $d$ and goes to $\infty$ when $d \to \infty$. Let us then consider the particular value of $d$ where this minimum is achieved. Solving the fixed point equation for this particular value, we have

$$r^* = \frac{4C^2}{n} \left( \sqrt{d} + \sqrt{d + 8\sqrt{n}R(4C)^{-1} \sqrt{\sum_{j > d} \lambda_j}} \right)^2.$$

Now for any other value $d' \neq d$, $r^*$ satisfies

$$r^* \leq \frac{4C}{\sqrt{n}} \left( \sqrt{r^* d'} + 2R \sqrt{\sum_{j > d'} \lambda_j} \right),$$

which means that $r^*$ is smaller than the largest solution of the corresponding equality. As a result, we have

$$r^* = \inf_{d \in \mathbb{N}} \frac{4C^2}{n} \left( \sqrt{d} + \sqrt{d + 8\sqrt{n}R(4C)^{-1} \sqrt{\sum_{j > d} \lambda_j}} \right)^2.$$

Using $(a + b)^2 \leq 2(a^2 + b^2)$, putting $C = C_R$ and finally using the assumption $RC_R^{-1} \leq \eta_1 M^{-1}$ yields the result.  $\square$

This concludes the proof of Theorem 4.4 for setting (S1). We finally turn to checking hypothesis (H4) in setting (S2): in this case we use a classical entropy chaining argument.

THEOREM 6.9. *Under assumption* (A1) *and the notation of setting* (S2), *we have*

$$\mathbb{E}\left[\sup_{\substack{\|g\|_k \leq R \\ d_2^2(g, g_0) \leq r}} |(P - P_n)(\ell(g) - \ell(g_0))|\right] \leq \frac{48R}{\sqrt{n}}\xi\left(\frac{\sqrt{r}}{2R}\right) + \frac{8MR^3}{n}r^{-1}\xi\left(\frac{\sqrt{r}}{2R}\right)^2$$

$$:= \psi_R(r),$$

*where the function* $\xi$ *is defined as in* (3.2). *The function* $\psi_R$ *is sub-root; if* $x_*$ *denotes the solution of the equation* $\xi(x) = M^{-1}n^{-1/2}x^2$, *then the solution* $r_R^*$ *of the equation* $\psi_R(r) = C_R^{-1}r$, *with* $C_R \geq MR$, *satisfies*

$$r_R^* \leq 2500M^{-2}C_R^2 x_*^2.$$

The chaining technique used for proving this theorem is summed up in the next lemma, for which we give a proof for completeness.

LEMMA 6.10. *Let* $\mathcal{F}$ *be a class of real functions which is separable in the supremum norm, containing the null function, and such that every* $f \in \mathcal{F}$ *satisfies* $\|f\|_\infty \leq M$ *and* $\mathbb{E}[f^2] \leq \sigma^2$. *Denote* $H_\infty(\varepsilon)$ *the supremum norm* $\varepsilon$-*entropy for* $\mathcal{F}$. *Then it holds that*

$$(6.19) \qquad \mathbb{E}\left[\sup_{f \in \mathcal{F}} |(P - P_n)f|\right] \leq \frac{24}{\sqrt{n}}\int_0^\sigma \sqrt{H_\infty(\varepsilon)}\, d\varepsilon + \frac{MH_\infty(\sigma)}{n}.$$

PROOF. It is a well-known consequence of Hoeffding's (resp. Bernstein's) inequality that a finite class of functions $\mathcal{G}$ bounded by $M$ in absolute value have

$$(6.20) \qquad \mathbb{E}\left[\sup_{g \in \mathcal{G}}(P - P_n)g\right] \leq \sqrt{2\frac{M^2 \log(|\mathcal{G}|)}{n}};$$

respectively, if, additionally, it holds that $\mathbb{E}[g^2] \leq \sigma^2$ for all $g \in \mathcal{G}$, we have

$$(6.21) \qquad \mathbb{E}\left[\sup_{g \in \mathcal{G}}(P - P_n)g\right] \leq \sqrt{2\frac{\sigma^2 \log(|\mathcal{G}|)}{n}} + \frac{M \log(|\mathcal{G}|)}{3n}.$$

Since $\mathcal{F}$ contains the null function, it is clear that

$$(6.22) \qquad \mathbb{E}\left[\sup_{f \in \mathcal{F}} |(P - P_n)f|\right] \leq \mathbb{E}\left[\sup_{f \in \mathcal{F}}(P - P_n)f\right] + \mathbb{E}\left[\sup_{f \in \mathcal{F}}(P_n - P)f\right].$$

Since we have assumed that $\mathcal{F}$ is separable for the sup norm, it is sufficient to prove (6.19) for any finite subset of $\mathcal{F}$. Without loss of generality, we therefore assume that $\mathcal{F}$ is finite. Put $\delta_i = \sigma 2^{-i}$ and let, for any $f \in \mathcal{F}$, $\Pi_i f$ be a member of a $\delta_i$-supremum norm cover of $\mathcal{F}$ such that $\|\Pi_i f - f\|_\infty \leq \delta_i$. We write

$$\mathbb{E}\Big[\sup_{f \in \mathcal{F}}(P - P_n)f\Big] \leq \mathbb{E}\Big[\sup_{f \in \mathcal{F}}(P - P_n)\Pi_0 f\Big]$$
$$+ \sum_{i>0}\mathbb{E}\Big[\sup_{f \in \mathcal{F}}(P - P_n)(\Pi_i f - \Pi_{i-1}f)\Big].$$

We now apply (6.21) to the first term of the above bound and (6.20) to all of the other terms. More precisely, we apply (6.21) to the class $\{\Pi_0 f, f \in \mathcal{F}\}$ which has cardinality bounded by $\exp(H_\infty(\delta_0))$, and respectively (6.20) to the classes $\{(\Pi_i f - \Pi_{i-1}f), f \in \mathcal{F}\}$ which have their respective cardinality bounded by $\exp(2H_\infty(\delta_i))$. We then have

$$\mathbb{E}\Big[\sup_{f \in \mathcal{F}}(P - P_n)f\Big] \leq \sqrt{\frac{2\sigma^2 H_\infty(\delta_0)}{n}} + \frac{MH_\infty(\delta_0)}{3n} + \sum_{i>0}\sqrt{\frac{36\delta_i^2}{n}H_\infty(\delta_i)}$$
$$\leq \frac{12}{\sqrt{n}}\int_0^\sigma \sqrt{H_\infty(\varepsilon)}\,d\varepsilon + \frac{MH_\infty(\sigma)}{3n}.$$

We apply the same inequality to the class $-\mathcal{F}$ and conclude using (6.22). $\square$

PROOF OF THEOREM 6.9.  We want to apply Lemma 6.10 to the class of functions

$$\mathcal{F}_{g_0} = \{\ell(g) - \ell(g_0) : g \in \mathcal{H}_k; \|g\|_k \leq R; \mathbb{E}[(\ell(g) - \ell(g_0))^2] \leq r\}.$$

Similarly to the reasoning used in the proof of Theorem 6.6, it is clear that

$$\mathcal{F}_{g_0} \subset \widetilde{\mathcal{F}}(2R, r) = \{\ell(g); g \in \mathcal{H}_k; \|g\|_k \leq 2R; \mathbb{E}[\ell(g)^2] \leq r\}.$$

Because the loss function $\ell$ is 1-Lipschitz, it holds that $\|\ell(f) - \ell(g)\| \leq \|f - g\|_\infty$, hence, $H_\infty(\widetilde{\mathcal{F}}(2R, r), \varepsilon) \leq H_\infty(\mathcal{B}_{\mathcal{H}_k}(2R), \varepsilon)$. Applying Lemma 6.10 therefore yields

$$\mathbb{E}\Big[\sup_{\substack{\|g\|_k \leq R \\ d_2^2(g, g_0) \leq r}} |(P - P_n)(\ell(g) - \ell(g_0))|\Big]$$
$$\leq \frac{24}{\sqrt{n}}\int_0^{\sqrt{r}} \sqrt{H_\infty(2R\mathcal{B}_{\mathcal{H}_k}, \varepsilon)}\,d\varepsilon + \frac{2MRH_\infty(2R\mathcal{B}_{\mathcal{H}_k}, \sqrt{r})}{n}.$$
$$= \frac{48R}{\sqrt{n}}\int_0^{\sqrt{r}/2R} \sqrt{H_\infty(\mathcal{B}_{\mathcal{H}_k}, \varepsilon)}\,d\varepsilon + \frac{2MRH_\infty(\mathcal{B}_{\mathcal{H}_k}, \sqrt{r}/2R)}{n}$$

$$\leq \frac{48R}{\sqrt{n}}\xi\left(\frac{\sqrt{r}}{2R}\right) + \frac{8MR^3}{n}r^{-1}\xi\left(\frac{\sqrt{r}}{2R}\right)^2 = \psi_R(r),$$

where $\xi$ is defined as in (3.2), and the last inequality comes from the observation that $\xi(x) \leq x\sqrt{H_\infty(\mathcal{B}_{\mathcal{H}_k}, x)}$. The function $\psi_R$ is obviously sub-root since $H_\infty(\mathcal{B}_{\mathcal{H}_k}, \varepsilon)$ is a decreasing function or $\varepsilon$.

Denote $x_*$ the solution of the equation $\xi(x) = M^{-1}\sqrt{n}x^2$; we claim that, for a suitable choice of constant $c$, $t_R^* = c^2 M^{-2} C_R^2 x_*^2$ is an upper bound for the solution $r_R^*$ of the equation $\psi_R(r) = C_R^{-1}r$ entering in hypothesis (H4). This is implied by the relation $\psi_R(t_R^*) \leq C_R^{-1}t_R^*$ which we now prove.

Note that $\frac{\sqrt{t_R^*}}{2R} = c\frac{C_R}{2RM}x_*$, and that $\frac{C_R}{RM} \geq 1$. Since $x^{-1}\xi(x)$ is a decreasing function, assuming $c \geq 2$, it holds that

$$\xi\left(\frac{\sqrt{t_R^*}}{2R}\right) \leq c\frac{C_R}{2RM}\xi(x_*) = c\frac{C_R\sqrt{n}}{2RM^2}x_*^2 = \frac{\sqrt{n}}{cRC_R}t_R^*.$$

Plugging this into the expression for $\psi_R$ yields

$$\psi_R(t_R^*) \leq \left(\frac{48}{c} + \frac{8MR}{c^2 C_R}\right)\frac{t_R^*}{C_R} \leq \left(\frac{48}{c} + \frac{8}{c^2}\right)\frac{t_R^*}{C_R},$$

where we have used again the relation $MR \leq C_R$. The choice $c = 50$ implies the desired relation. $\square$

6.4. *Proof of Theorem* 3.1. Theorem 4.4 states that the conditions (H1)–(H4) of the model selection theorem (Theorem 4.3) are satisfied for the family of models $\mathcal{B}(R), R \in \mathcal{R}$ and some explicit values for $b_R, C_R, \phi_R$ and $r_R^*$ [depending on the considered setting (S1) or (S2)]. Let us choose an appropriate finite set $\mathcal{R}$ and a sequence $(x_R)_{R \in \mathcal{R}}$ so that we can approximate the minimization over all $R > 0$ in equation (3.4) by a minimization over the finite set of radii $\mathcal{R}$.

We consider the set of discretized radii

$$\mathcal{R} = \{M^{-1}2^k, k \in \mathbb{N}, 0 \leq k \leq \lceil \log_2 n \rceil\}.$$

The cardinality of $\mathcal{R}$ is then $\lceil (\log_2 n) \rceil + 1$ and we consequently choose $x_R \equiv \log(\log_2 n + 2)$ for all $R \in \mathcal{R}$ which satisfies $\sum_{R \in \mathcal{R}} e^{-x_R} \leq 1$.

In order to apply Theorem 4.3, the penalty function should satisfy equation (4.2). A sufficient condition on the penalty function for the family of models $\{\mathcal{B}(R), R \in \mathcal{R}\}$ is therefore

$$\text{pen}(R) \geq c_1\left(\frac{r_R^*}{C_R} + \frac{(C_R + b_R)(x_R + \log(\delta^{-1}) \vee 1)}{n}\right),$$

where $c_1$ is a suitable constant, and we picked $K = 3$ in equation (4.2).

Recalling the definition of $\gamma(n)$ in settings (S1) and (S2), the requirement (3.3) on $\Lambda_n$ and the definition of $w_1$ in Theorem 3.1, it can be checked by elementary manipulations that the above condition on the penalty is satisfied in both settings for

$$\text{pen}(R) = \Lambda_n(\varphi(MR/2) + w_1\eta_0^{-1}),$$

up to a suitable choice of the constant $c$ in (3.3); note that we can assume $c \geq 1$.

The last step to be analyzed now is how to go back and forth between the discretized framework $R \in \mathcal{R}$ and the continuous framework to obtain the final result. To apply the model selection theorem, we will interpret the continuous regularization defining $\widehat{g}$ as an approximate discretized penalized minimization over the above family of models using the penalty function defined above.

In view of definition (3.4) of the estimator $\widehat{g}$, the following upper bound holds:

$$P_n\ell(\widehat{g}) + \Lambda_n\varphi(M\|\widehat{g}\|_k) \leq P_n\ell(0) + \Lambda_n\varphi(0) = 1,$$

which implies $1 \geq \Lambda_n\varphi(M\|\widehat{g}\|_k)$. Since we have assumed $c \geq 1$ in (3.3), we have $\Lambda_n \geq n^{-1}$; this implies $\|\widehat{g}\|_k \leq M^{-1}n$ (using the assumption on $\varphi$). Denote $\widehat{R} = M^{-1}2^{\widehat{k}}$ where $\widehat{k} = \lceil(\log_2(M\|\widehat{g}\|_k))_+\rceil$. The fact that $\|\widehat{g}\|_k \leq M^{-1}n$ implies $\widehat{R} \in \mathcal{R}$. Note that $\widehat{g} \in \mathcal{B}(\widehat{R})$ and that $\widehat{R} \leq 2M^{-1}\max(M\|\widehat{g}\|_k, 1)$. This entails

$$\begin{aligned}
P_n\ell(\widehat{g}) + \text{pen}(\widehat{R}) &\leq P_n\ell(\widehat{g}) + \Lambda_n\varphi(\max(M^{-1}\|\widehat{g}\|_k, 1)) + \eta_0^{-1}w_1^{-1}\Lambda_n \\
&\leq P_n\ell(\widehat{g}) + \Lambda_n\varphi(M^{-1}\|\widehat{g}\|_k) + \Lambda_n\varphi(1) + w_1^{-1}\eta_0^{-1}\Lambda_n \\
&\leq \inf_{g \in \mathcal{H}_k}[P_n\ell(g) + \Lambda_n\varphi(M^{-1}\|g\|_k)] + \Lambda_n\varphi(1) + w_1^{-1}\eta_0^{-1}\Lambda_n \\
&= \inf_{R \geq 0}\inf_{g \in \mathcal{B}(R)}[P_n\ell(g) + \Lambda_n(\varphi(MR) + \varphi(1) + w_1^{-1}\eta_0^{-1})] \\
&\leq \inf_{R \in \mathcal{R}}\inf_{g \in \mathcal{B}(R)}[P_n\ell(g) + \Lambda_n(\varphi(MR) + \varphi(1) + w_1^{-1}\eta_0^{-1})],
\end{aligned}$$

where the first inequality follows from the definition of $\text{pen}(\widehat{R})$, and the third from the definition of $\widehat{g}$. So if we put $\rho_R = \Lambda_n(\varphi(MR) + \varphi(1) + w_1^{-1}\eta_0^{-1}) - \text{pen}(R) \geq 0$, we just proved that $\widehat{g}$ is a $(\rho_R)$-approximate penalized minimum loss estimator over the family $(\mathcal{B}(R))_{R \in \mathcal{R}}$. Now applying the model selection theorem (Theorem 4.3), we conclude that the following bound holds with probability at least $1 - \delta$:

$$\begin{aligned}
L(\widehat{g}, s^*) &\leq 2\inf_{R \in \mathcal{R}}\inf_{g \in \mathcal{B}(R)}[L(g, s^*) + 2\Lambda_n(\varphi(MR) + \varphi(1) + w_1^{-1}\eta_0^{-1})] \\
&= 2\inf_{R \in \mathcal{R}}\inf_{g \in \mathcal{B}(R)}[L(g, s^*) + 2\Lambda_n\varphi(2^{\log_2 MR})] + 4\Lambda_n(\varphi(1) + w_1^{-1}\eta_0^{-1}) \\
&\leq 2\inf_{M^{-1} \leq R \leq nM^{-1}}\inf_{g \in \mathcal{B}(R)}[L(g, s^*) + 2\Lambda_n\varphi(2^{\lceil\log MR\rceil})]
\end{aligned}$$

$$+ 4\Lambda_n(\varphi(1) + w_1^{-1}\eta_0^{-1})$$

$$\leq 2 \inf_{R \leq nM^{-1}} \inf_{g \in \mathcal{B}(R)} [L(g, s^*) + 2\Lambda_n \varphi(2^{\lceil (\log MR)_+ \rceil})]$$

$$+ 4\Lambda_n(\varphi(1) + w_1^{-1}\eta_0^{-1})$$

$$\leq 2 \inf_{R \leq nM^{-1}} \inf_{g \in \mathcal{B}(R)} [L(g, s^*) + 2\Lambda_n(\varphi(2(MR \vee 1)))]$$

$$+ 4\Lambda_n(\varphi(1) + w_1^{-1}\eta_0^{-1})$$

$$\leq 2 \inf_{g \in \mathcal{H}_k} [L(g, s^*) + 2\Lambda_n \varphi(2M\|g\|_k)] + 4\Lambda_n(2\varphi(2) + w_1^{-1}\eta_0^{-1}).$$

The last inequality holds because, if we denote $g^*$ the minimizer of the last infimum, comparing it with the constant null function (as for $\hat{g}$ earlier), we conclude that $2\Lambda_n \varphi(2M\|g^*\|_k) \leq 1$, implying $\|g^*\|_k \leq M^{-1}n$, so that the restriction $R \leq M^{-1}n$ in the previous infimum can be dropped.

## APPENDIX A: PROPERTIES OF THE KERNEL INTEGRAL OPERATOR.

In this appendix, we sum up a few useful properties of the integral operator $L_k$ introduced in (3.1). These are used in the proof of Lemma 6.8. While these results are certainly not new, we provide a self-contained proof for completeness.

LEMMA A.1. *Let $\mathcal{H}_k$ be a separable RKHS with kernel $k$ on a measurable space $\mathcal{X}$. Assume $y \mapsto k(x, y)$ is measurable for any fixed $x \in \mathcal{X}$. Then the function $x \mapsto k(x, \cdot) \in \mathcal{H}_k$ is measurable; in particular, $(x, y) \mapsto k(x, y)$ is jointly measurable.*

*Let $P$ be a probability distribution on $\mathcal{X}$; assume $L^2(P)$ is separable and $\mathbb{E}_{X \sim P}[k(X, X)] < \infty$.*

*Then $\mathcal{H}_k \subset L^2(P)$ and the canonical inclusion $T : \mathcal{H}_k \to L^2(P)$ is continuous.*

*The integral operator $L_k : L^2(P) \to L^2(P)$ defined as*

$$(L_k f)(x) = \int k(x, y) f(y) \, dP(y)$$

*is well defined, positive, self-adjoint and trace class; moreover, $L_k = TT^*$. In particular, if $(\lambda_i)_{i \geq 0}$ denote its eigenvalues, repeated with their multiplicities, $\sum_{i \geq 0} \lambda_i < \infty$.*

*Finally, there exists an orthonormal basis $(\psi_i)_{i \geq 0}$ of $\mathcal{H}_k$ such that, for any $f \in \mathcal{H}_k$,*

$$\|Tf\|_{2,P}^2 = \sum_{i \geq 0} \lambda_i \langle f, \psi_i \rangle^2.$$

PROOF. Let us first prove that any function $f \in \mathcal{H}_k$ is measurable. By assumption, for any fixed $x$, $k(x, \cdot)$ is measurable; hence, also any finite linear combination $\sum_i \alpha_i k(x_i, \cdot)$. Any function $f \in \mathcal{H}_k$ is the limit in $\mathcal{H}_k$ of a sequence of such linear combinations. By the reproducing property, a sequence converging in $\mathcal{H}_k$ also converges pointwise, since $\langle f_i, k(x, \cdot) \rangle = f_i(x)$. Hence, $f$ is measurable. Now we prove that $x \mapsto K(x) = k(x, \cdot) \in \mathcal{H}_k$ is measurable.

For any $f \in \mathcal{H}_k$, $x \mapsto \langle k(x, \cdot), f \rangle = f(x)$ is measurable, hence, the inverse image of a half-space by $K$ is measurable. Since $\mathcal{H}_k$ is separable, any open set is a countable union of open balls (Lindelöf property); and any ball in $\mathcal{H}_k$ is a countable intersection of half-spaces. Hence, $K$ is measurable. This implies that $k(x, y) = \langle k(x, \cdot), k(y, \cdot) \rangle_{\mathcal{H}}$ is jointly measurable.

By the Cauchy–Schwarz inequality, we further have $|k(x, y)|^2 \leq k(x, x)k(y, y)$, so that the assumptions that $k(x, x) \in L_1(P)$ imply that $k(\cdot, \cdot) \in L^2(\mathcal{X} \times \mathcal{X}, P \otimes P)$. This ensures that $L_k$ is well defined [as an operator $L_P^2(\mathcal{X}) \to L_P^2(\mathcal{X})$] and Hilbert–Schmidt, hence, compact. Moreover, by symmetry of $k$, $L_k$ is self-adjoint. As $L^2(P)$ is separable, $L_k$ can be diagonalized in an orthonormal basis $(\phi_i)_{i \geq 0}$ of $L^2(P)$, where $L_k \phi_i = \lambda_i \phi_i$.

Consider now the canonical inclusion $T$ from the reproducing kernel Hilbert space $\mathcal{H}_k$ into $L^2(P)$. For $f \in \mathcal{H}_k$, we have

$$\int f^2(x)\, dP(x) = \int \langle f, k(x, \cdot) \rangle_{\mathcal{H}_k}^2\, dP(x) \leq \|f\|_{\mathcal{H}_k}^2 \int k(x, x)\, dP(x).$$

This proves that $T$ is well defined and continuous on $\mathcal{H}_k$. Let $T^* : L^2(P) \to \mathcal{H}_k$ denote its adjoint.

For any $f \in L^2(P)$, we have by definition for all $x \in \mathcal{X}$, $T^* f(x) = \langle k(x, \cdot), T^* f \rangle_{\mathcal{H}_k} = \langle Tk(x, \cdot), f \rangle_{L^2(P)} = (L_k f)(x)$. Hence, $TT^* = L_k$. In particular, $\lambda_i = \langle \phi_i, \lambda_i \phi_i \rangle_{L^2(P)} = \langle T^* \phi_i, T^* \phi_i \rangle_{\mathcal{H}_k} \geq 0$, which proves that $L_k$ is a positive operator.

Now let us consider the operator $C = T^*T : \mathcal{H}_k \to \mathcal{H}_k$. It is bounded, positive and self-adjoint. Let $(\psi_i)_{i \geq 0}$ be an orthonormal basis of $\mathcal{H}_k$. We have

$$k(x, x) = \langle k(x, \cdot), k(x, \cdot) \rangle = \sum_{i \geq 0} \langle k(x, \cdot), \psi_i \rangle^2 = \sum_{i \geq 0} \psi_i(x)^2$$

and

$$\sum_{i \geq 0} \langle \psi_i, C\psi_i \rangle_{\mathcal{H}_k} = \sum_{i \geq 0} \|T\psi_i\|_{2,P}^2 = \mathbb{E} \sum_{i \geq 0} (T\psi_i)^2(X) = \mathbb{E} k(X, X) < \infty,$$

by monotone convergence. This proves that $C$ is trace-class. Now since $TT^*$ and $T^*T$ have the same nonzero eigenvalues (with identical multiplicities), and $\mathrm{tr}\, C = \sum_{i \geq 0} \lambda_i < \infty$, $L_k$ is also trace-class.

We can actually choose $(\psi_i)$ as an orthonormal basis of eigenvectors of $C$ with corresponding eigenvalues $\lambda_i$. In that case, we can write any function $f \in \mathcal{H}_k$ as

$$f = \sum_{i \geq 0} \langle f, \psi_i \rangle \psi_i,$$

where $\|f\|^2_{\mathcal{H}_k} = \sum_{i \geq 0} \langle f, \psi_i \rangle^2$ and by continuity of $T$,

$$Tf = \sum_{i \geq 0} \langle f, \psi_i \rangle T\psi_i.$$

Now, since $C\psi_i = \lambda_i \psi_i$, we have $\langle T\psi_i, T\psi_j \rangle = \lambda_i \langle \psi_i, \psi_j \rangle = \lambda_i \delta_{ij}$ so that

$$\|Tf\|^2_{2,P} = \sum_{i \geq 0} \lambda_i \langle f, \psi_i \rangle^2. \qquad \square$$

**Acknowledgments.** We are grateful to Jean–Yves Audibert and Nicolas Vayatis for helpful comments and pointing out some mistakes in earlier versions of this work. We thank the two anonymous reviewers who gave many insightful comments and helped us improving our initial results and the overall organization of the paper.

G. Blanchard
Fraunhofer–Insitute FIRST
Kekuléstr. 7
12489 Berlin
Germany
E-mail: blanchar@first.fhg.de

O. Bousquet
Google
8002 Zurich
Switzerland
E-mail: olivier.bousquet@gmail.com

P. Massart
Université Paris-Sud
Laboratoire de Mathématiques
91405 Orsay Cedex
France
E-mail: massart@stats.math.u-psud.fr